\documentclass{amsart}
\usepackage{a4wide}

\usepackage{amssymb}
\usepackage{amsmath}
\usepackage[utf8]{inputenc}
\usepackage{graphicx}

\usepackage{color}
\usepackage{siunitx}
\usepackage{comment}
\usepackage{hyperref}
\usepackage{tikz}
\usetikzlibrary{calc}
\def\centerarc[#1](#2)(#3:#4:#5)
{ \draw[#1] ($(#2)+({#5*cos(#3)},{#5*sin(#3)})$) arc (#3:#4:#5); }

\newcommand{\tsr}[1]{{\boldsymbol{\mathrm{#1}}}}
\newcommand{\eps}{\varepsilon}
\renewcommand{\vec}{\boldsymbol}
\newcommand{\opdiv}{\operatorname{div}}

\newcommand{\thickness}{t_\Surf}

\newcommand{\Surf}{\mathcal{S}}
\newcommand{\Bulk}{\mathcal{B}}

\newcommand{\Cof}{\operatorname{Cof}}
\newcommand{\trace}{\operatorname{tr}}

\newcommand{\hatkappa}{\hat{\tsr \kappa}}
\newcommand{\hatm}{\hat{\tsr m}}

\title{Direct coupling of continuum and shell elements in large deformation problems}

\author{Astrid Pechstein} 
\address{Astrid Pechstein, Institute of Technical Mechanics, Johannes Kepler University Linz, Altenberger Str. 69, 4040 Linz, Austria}
\email{astrid.pechstein@jku.at}

\author{Michael Neunteufel} 
\address{Michael Neunteufel, Fariborz Maseeh Department of Mathematics and Statistics, Portland State University, 1825 SW Broadway, 97201 Portland, Oregon, USA}
\email{mneunteu@pdx.edu}
           

\begin{document}

\begin{abstract}
In many applications, thin shell-like structures are integrated within or attached to volumetric bodies. This includes reinforcements placed in soft matrix material in lightweight structure design, or hollow structures that are partially or completely filled.
Finite element simulations of such setups are highly challenging. A brute force discretization of structural as well as volumetric parts using well-shaped three-dimensional elements may be accurate, but leads to problems of enormous computational complexity even for simple models. One desired alternative is the use of shell elements for thin-walled parts, as such a discretization greatly alleviates size restrictions on the underlying finite element mesh. However, the coupling of different formulations within a single framework is often not straightforward and may lead to locking if not done carefully.
Neunteufel and Schöberl proposed a mixed shell element where, apart from displacements of the center surface, bending moments are used as independent unknowns. These elements were not only shown to be locking free and highly accurate in large-deformation regime, but also do not require differentiability of the shell surface.
They can directly be coupled to classical volume elements of arbitrary order by sharing displacement degrees of freedom at the center surface, thus achieving the desired coupled discretization. As the elements can be used on unstructured meshes, adaptive mesh refinement based on local stress and bending moments can be used. We present computational results that confirm exceptional accuracy for problems where thin-walled structures are embedded as reinforcements within soft matrix material.   
\end{abstract}

\maketitle



\section{Introduction}

As lightweight structures are designed, one often meets with stiff thin-walled structures (the \emph{reinforcement}) in combination with a soft matrix. Reinforcements may be submerged within the matrix, as well as form the surface of the lightweight structure \cite{CBG2020,KK2005,KS2019}. Finite element simulations of these combined structures prove to be computationally expensive due to their adverse geometric dimensions. As thin-walled structures are resolved by conventional elements, the thickness of the reinforcements makes small mesh sizes in these regions necessary, which leads to an excessive growth in the number of degrees of freedom.

To alleviate the restrictions on the mesh size, solid shell elements have been developed especially for the discretization of thin structures. Optimally, the reinforcement can then be discretized using only one layer of solid shell elements, where the mesh size in thickness direction is independent of the in-plane size of these elements. This greatly reduces the number of elements, and thus the number of unknowns to be computed. One well-known approach is the use of assumed natural strains as additional, element-internal degrees of freedom in the solid-shell elements, as was proposed by Sze and Yao \cite{Sze:2000}, and extended to the large deformation case in \cite{Sze:2002}. Another large deformation element based on the Hu--Washizu three-field principle is presented by Klinkel et al.~\cite{Klinkel:2006}. 

As these elements are constructed in such a way that degrees of freedoms and shape functions coincide with conventional elements on the reinforcement/matrix interface, coupling of the different element types in a combined structure is straightforward when using a matching pair of solid shell and standard elements. Limitations of this approach are found at the boundaries of flat reinforcement structures: here, the thickness still limits the mesh size (although much less so than in the case of standard elements), and the connection of several shells along a common line (e.g. a T-joint) is not straightforward, especially if the shells do not meet at an angle of $\SI{90}{\degree}$.

A different approach is the modeling of thin-walled reinforcements as shells, which are included as surface or interface within the soft matrix. This approach shows great advantages from the point of view of mesh generation, as then the thickness does not appear in the finite element mesh at all. Reinforcements are modelled as shell, only its mid-surface is discretized.
Under the kinematic assumptions of Kirchhoff--Love theory, elastic energy on the shell mid-surface can be derived analytically, leading to fourth order partial differential equations for the displacement, see e.g. the monograph \cite{Vetyukov:2014book}. Finite elements have to satisfy $\mathcal{C}^1$ continuity, where no kinks in the displacement are allowed. For the geometrically linear case, one well-known example is the Argyris triangle from the TUBA family \cite{Argyris:1968}. Extensions to large deformations were, however, of limited success \cite{Argyris:1969}. A quadrilateral element suitable for large deformations was proposed by Vetyukov \cite{Vetyukov:2014}. It is an element from the family of absolute nodal coordinate formulation based elements, and uses positions and assorted slope vectors as degrees of freedom. Additionally, we cite the triangular six-node element by Dau et al.~\cite{Dau:2006} and the monograph by Bernadou and James~\cite{Bernadou:1996} for different Kirchhoff--Love shell element formulations. Common to all these specific shell elements is the fact that the degrees of freedom for the mid-surface displacement do not match the nodal degrees of freedom of elements used in the matrix. Thus, coupling of these discretizations is not straightforward, and leads to locking or loss of accuracy if not done carefully.




In the field of isogeometric analysis (IGA), splines of higher smoothness are used as ansatz functions in the discretization, see the monograph by Cottrell, Hughes, and Bazilevs~\cite{Cottrell:2009} for an introduction. Splines lend themselves naturally to the discretization of Kirchhoff--Love shells, as the $\mathcal{C}^1$ differentiability can be ensured through the choice of order and multiplicity of knots in the discretization. Bazilevs et al.~\cite{Bazilevs:2010} proposed to use T-splines for volumes as well as plate domains, for solids as well as fluid mechanics. Interpolation by non-uniform rational B-splines (NURBS)  was proven to be effective for modeling classical shells by Kiendl et al.~\cite{Kiendl:2009}, and extended to compressible and incompressible hyperelastic materials in \cite{Kiendl:2015}. Regarding reinforcements, shell elements for fibre-reinforced composites were designed by Schulte and co-workers~\cite{Schulte:2020}. When the same class of IGA elements can be used for  reinforcement and matrix, a direct coupling of the approaches is possible, see Liu et al.~\cite{Liu:2021}. However, the discretization of complex reinforced structures requires a discretization consisting of multiple patches. For these multiple patches, $\mathcal{C}^1$ continuity is lost, and needs to be re-established through additional constraints. Schuß et al.~\cite{Schuss:2019} presented a methodology applicable to Kirchhoff--Love shells. However, the joint discretization of matrix volume and shell structure using IGA splines is certainly not straightforward. Also, the treatment of kinks in the surface or joints of more than one shell along a common line is to be done with care. Last, correct representations of the coupling of trimmed surfaces and volumes has not been treated in the literature to the best knowledge of the authors.

Pechstein and Schöberl introduced the \emph{Tangential-Displacement-Normal-Normal-Stress} (TDNNS) method \cite{Pechstein:2011} and proved that it is locking-free in the sense of solid-shell elements for the small deformation case in \cite{Pechstein:2012}. An extension to the large deformation case could be established by Neunteufel et al.~\cite{Neunteufel:2021threefield} by way of a three-field Hu--Washizu formulation. Concerning structural elements, a plate formulation for Reissner--Mindlin plates based on these techniques was presented in \cite{Pechstein:2017}. These ideas were generalized to large-deformation shells by Neunteufel and Schöberl in ~\cite{Neunteufel:2019}, where a family of arbitrary order mixed shell elements was proposed. When reducing this shell formulation to linear Kirchhoff--Love plates, the mixed Hellan--Herrmann--Johnson method is recovered \cite{Hel67,Her67, Joh73, Com89}. Membrane locking of the original method \cite{Neunteufel:2019} was treated efficiently in \cite{Neunteufel:2021}. Coupling of these shell elements with electromechanics for dielectric elastomers has been successfully performed by Pechstein and Krommer in \cite{PK2024}. Essential for all these formulations based on the original contribution \cite{Neunteufel:2019} is that in this formulation $\mathcal{C}^1$ continuity for the displacement could be alleviated by the introduction of an independent moment tensor. Nodal elements without additional smoothness can directly be used to discretize the shell mid-surface displacement. Thus, a coupled discretization of reinforcements and matrix is possible in a straightforward manner.

The coupling of shear deformable Reissner--Mindlin shells to solid domains has been treated in \cite{YYM2019,KF2024}. Coupling of Cosserat shells in the framework of tangential calculus is discussed by Sky et al.~in \cite{SHZBN2024}, where additional rotational degrees of freedom are used and needed to be coupled with solid (Cosserat) elasticity. The additional degrees of freedom in these shell formulations alleviate the mid-surface smoothness constraint, but can lead to (shear) locking, which has to be taken into consideration. The extension to shear deformable shell models is, however, out of scope of this work.

The present paper is organized as follows: In the next section we will describe the original elasticity problem involving structure and bulk domains. Section~\ref{sec:geometric_kinematic_reduction} is devoted to the geometric and kinematic reduction of the full elasticity problem corresponding shell and a bulk part. Therein the continuous shell equations are derived. In Section~\ref{sec:mixed_shell_formulation} we present a consistent discretization by mixed shell elements. We discuss involved finite element spaces, hybridization, and the applicability to structures with kinks and branches such as T-joints. Local equilibrium conditions are derived in Section~\ref{sec:local_equilibrium} for the case of small deformations and plates. In Section~\ref{sec:computational_results} we present several numerical examples to demonstrate the excellent performance and efficiency of the proposed method.





\section{Definition of the original problem}

\paragraph{Geometry and kinematics of the original problem}

Let $\Omega \subset \mathbb R^3$ denote the domain associated to the so-called reference configuration of the body of interest. In the following, we assume that this domain can be subdivided into two parts: the shell domain, which is of small thickness, but otherwise arbitrary topology concerning its mid-surface, and the bulk domain, which denotes the remaining part of $\Omega$. To be concise, let $\Surf \subset \mathbb R^3$ be an oriented surface, i.e.~a two-dimensional manifold with normal vector $\vec N_\Surf$. Then, for given thickness $\thickness$, we define the shell domain $\Omega_\Surf$ as
\begin{align}
    \Omega_\Surf &= \{\vec X + \zeta \vec N_\Surf: \vec X \in \Surf, \zeta \in [-\thickness/2, \thickness/2]\}.
\end{align}
It is assumed that the mid-surface $\Surf$ and thickness $\thickness$ correspond insofar as that $\Omega_\Surf \subset \overline \Omega$; note that for the current definition, the shell domain may be immersed within $\Omega$, extend to the boundary $\partial \Omega$, and may also coincide with the original domain's surface. In all cases, the bulk domain is defined as the remaining part of $\Omega$,
\begin{align}
    \Omega_\Bulk &= \Omega \backslash \Omega_\Surf.
\end{align}

Let now $\vec u: \Omega \to \mathbb R^3$ denote the displacement vector describing the deformation of the body of interest. Using standard notation, we introduce $\nabla = \partial/\partial \vec X$ as the derivative with respect to material coordinates $\vec X$, and $\tsr F = (\nabla \vec u)^T + \tsr I$ the deformation gradient. Its determinant is denoted as $J = \det \tsr F$, while right Cauchy-Green tensor $\tsr C = \tsr F^T \cdot \tsr F$ and Green-Lagrange strain $\tsr E = 1/2(\tsr C - \tsr I)$ describe local deformations of the body. Last, any surface normal vector is mapped through the contravariant transformation,
\begin{align}
    \vec n_\Surf &= \frac{\Cof \tsr F \cdot \vec N_\Surf}{|\Cof \tsr F \cdot \vec N_\Surf|},
    \label{eq:transformN}
\end{align}
with $\Cof \tsr F = J \tsr F^{-T}$ the cofactor matrix. Above and in the following, we adhere to the convention to use upper-case letters for quantities associated to reference configuration, and lower-case letters for the corresponding quantities in deformed configuration.

\paragraph{Hyperelastic materials}
We assume that the elastic behavior of the solid is defined in the way of elastic energy densities. Let $\psi(\tsr C)$ denote the objective elastic energy density for the unreduced, original body. For presentation of the proposed method, we assume that all external forces are conservative body loads, such that there exists an energy density representation $\psi_{\mathrm{ext}}(\vec u)$. Of course, the generalization to surface loads as well as non-conservative loads is straightforward. Then, the static solution $\vec u$ satisfies the minimization problem
\begin{align}
    \int_\Omega \left(\psi(\tsr C) - \psi_\mathrm{ext}(\vec u) \right) \, dV \to \min_{\vec u \text{ adm.}}.
\end{align}
Above, any displacement field is \emph{admissible} if it corresponds to a finite elastic energy, $\int_\Omega \psi(\tsr C)\, dV < \infty$, and the kinematic boundary conditions on $\vec u$ are satisfied. In the context of linear elasticity, these requests usually result in the common choice $\vec u \in [H^1(\Omega)]^3$ with $\vec u = \vec 0$ on the Dirichlet boundary $\Gamma_\text{fix} \subset \partial \Omega$.
The principle of virtual work yields a variational formulation; let the solution $\vec u$ satisfy
\begin{align}
    \delta W_\text{int} = \int_\Omega \left( \frac{\partial \psi}{\partial \tsr C} : \delta \tsr C\right)\, dV &= 
    \int_\Omega \left( \frac{\partial \psi_\mathrm{ext}}{\partial \vec u} \cdot \delta \vec u \right)\, dV = \delta W_\mathrm{ext},
    \label{eq:PVA}
\end{align}
for \emph{admissible} virtual displacements $\delta \vec u$ satisfying homogeneous kinematic boundary conditions. Last, we state that the second and first Piola--Kirchhoff stress tensors are defined through the hyperelastic potential,
\begin{align}
  \tsr S &= 2\frac{\partial \psi}{\partial \tsr C}, &  \tsr P &= \frac{\partial \psi}{\partial \tsr F} = \tsr F \cdot \tsr S.
\end{align}

\paragraph{The case of small deformations}
If all deformations are small, reference and spatial configuration of the body of interest are commonly identified, and the Green strain tensor $\tsr E$ is replaced by its linearization
$\tsr \eps = 1/2(\nabla \vec u + (\nabla \vec u)^T)$. Cauchy stress $\tsr \sigma$ and Piola--Kirchhoff stresses $\tsr S$, $\tsr P$ are no longer distinguished, and a linear relationship between stress tensor $\tsr \sigma$ and strain $\tsr \eps$ is assumed, such that with the fourth order elasticity tensor $\mathbb C$,
\begin{align}
  \psi &= \frac{1}{2} \tsr \eps : \mathbb C : \tsr \eps, &
  \tsr \sigma &= \frac{\partial \psi}{\partial \tsr \eps} = \mathbb C : \tsr \eps.
  \label{eq:linearstress}
\end{align}

\section{Geometric and kinematic reduction}
\label{sec:geometric_kinematic_reduction}
\begin{figure}[ht]
  \centering

\begin{tikzpicture}[scale=0.525]
  \def\L{10} 
  \def\H{3} 
  \def\tS{0.6} 

  \def\offset{1.5*\L}
  
  \draw[fill=gray!30] (0,0) -- ++(\L,0) -- ++(0,\H) -- ++(-\L,0) -- cycle;
  \draw[fill=red!80] (0,{0.5*\H-0.5*\tS}) -- ++(\L,0) -- ++(0,\tS) -- ++(-\L,0) -- cycle;
      
  \draw[<->] (-0.5,0) -- (-0.5,\H) node[midway, left] {$L$};
  \draw[<->] (\L*1.05,{0.5*\H-0.5*\tS}) -- (\L*1.05,{0.5*\H+0.5*\tS}) node[midway, right] {$t$};

  \draw (0.3*\L,0.8*\H) node {$\Omega_\Bulk$};
  \draw (0.3*\L,0.2*\H) node {$\Omega_\Bulk$};
  \draw[color=red] (0.6*\L,{0.5*\H-1.5*\tS}) node {$\Omega_\Surf$};
 
  \draw[->] (\L+0.7,0.8*\H) to [out=60,in=120] (\offset-0.7,0.8*\H);

  \draw[fill=gray!30] (\offset,0) -- ++(\L,0) -- ++(0,\H) -- ++(-\L,0) -- cycle;
  \draw[color=red,very thick] (\offset,{0.5*\H}) -- (\offset+\L,{0.5*\H});
      
  \draw[<->] (\offset-0.5,0) -- (\offset-0.5,\H) node[midway, left] {$L$};

  \draw (\offset+0.3*\L,0.8*\H) node {$\Bulk$};
  \draw (\offset+0.3*\L,0.2*\H) node {$\Bulk$};
  \draw[color=red] (\offset+0.6*\L,{0.5*\H-1.5*\tS}) node {$\Surf$};
  
  \def\offsety{-1.6*\H}
  \def\tS{0.6} 

  \def\offset{1.5*\L}
  
  \draw[fill=gray!30] (0,\offsety) -- ++(\L,0) -- ++(0,\H) -- ++(-\L,0) -- cycle;
  \draw[fill=red!80] (0,{-0.5*\tS+\offsety}) -- ++(\L,0) -- ++(0,\tS) -- ++(-\L,0) -- cycle;
      
  \draw[<->] (-0.5,{-0.5*\tS+\offsety}) -- (-0.5,\H+\offsety) node[midway, left] {$L$};
  \draw[<->] (\L*1.05,{-0.5*\tS+\offsety}) -- (\L*1.05,{0.5*\tS+\offsety}) node[midway, right] {$t$};
  
  \draw (0.4*\L,{0.63*\H+\offsety}) node {$\Omega_\Bulk$};
  \draw[color=red] (0.6*\L,{-2*\tS+\offsety}) node {$\Omega_\Surf$};
  
  \draw[->] (\L+0.7,{0.8*\H+\offsety}) to [out=60,in=120] (\offset-0.9,0.8*\H+\offsety);
  
  \draw[fill=gray!30] (\offset,\offsety) -- ++(\L,0) -- ++(0,\H) -- ++(-\L,0) -- cycle;
  \draw[color=red,very thick] (\offset,\offsety) -- (\offset+\L,\offsety);
  
  \draw[<->] (\offset-0.5,{-0.5*\tS+\offsety}) -- (\offset-0.5,\H+\offsety) node[midway, left] {$L$};
  \draw[<-] (\offset+\L*1.05,{-0.5*\tS+\offsety}) -- (\offset+\L*1.05,{-1.5*\tS+\offsety}); 
  \draw[|-|] (\offset+\L*1.05,{-0.5*\tS+\offsety}) -- (\offset+\L*1.05,{\offsety}) node[midway, right] {$t/2$};
  \draw[->] (\offset+\L*1.05,{1*\tS+\offsety}) -- (\offset+\L*1.05,\offsety); 

  \draw (\offset+0.4*\L,0.63*\H+\offsety) node {$\Bulk$};
  \draw[color=red] (\offset+0.6*\L,{-1.5*\tS+\offsety}) node {$\Surf$};
\end{tikzpicture}
  \caption{Geometric reduction of the shell domain $\Omega_\Surf$ to a representation through its mid-surface~$\Surf$. Top: The shell domain is embedded into the bulk domain. Bottom: The shell domain is located at the border of the bulk domain.}
  \label{fig:geometricreduction}
\end{figure}
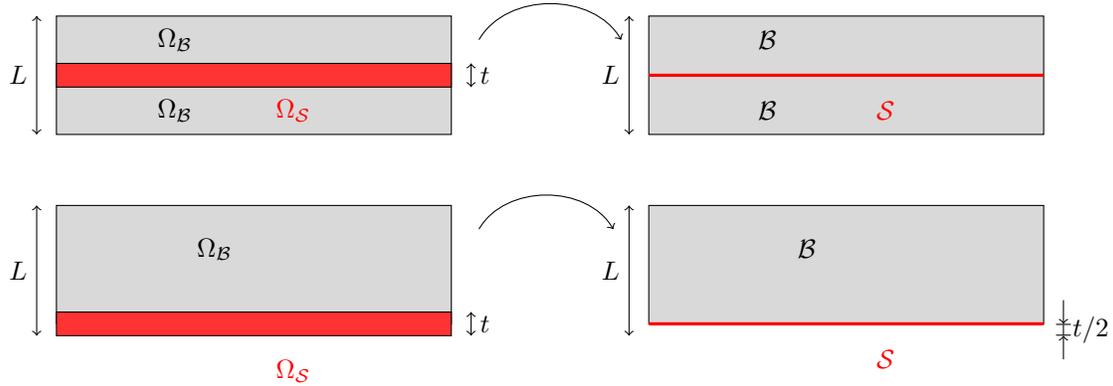

\paragraph{Geometric reduction}
The aim of the following considerations is to find a simplified approximate representation of the body's geometry that allows for a dimensional reduction of the shell contributions. The main idea is to reduce the shell domain $\Omega_\Surf$ to be represented through its mid-surface $\Surf$. In doing so, gaps are introduced in between the shell part's mid-surface $\Surf$ and the bulk domain $\Omega_\Bulk$. To restore the computational domain to be connected, the bulk domain $\Omega_\Bulk$ is modified such that it contains those parts of $\Omega_\Surf$ up to $\Surf$ as an interface or boundary. We introduce $\Bulk$ as this extended bulk domain; for a sketch see Figure~\ref{fig:geometricreduction}. Note that, if the shell forms part of the original geometry's boundary, the bulk $\Bulk$ covers those parts of $\Omega_\Surf$ \emph{inside} $\Surf$, while those parts \emph{outside} $\Surf$ are not included, and will not be discretized in simulations.


Thus, we are finally concerned with finding a displacement field $\vec u_\Bulk: \Bulk \to \mathbb R^3$ and $\vec u_\Surf: \Surf \to \mathbb R^3$ for the \emph{approximate} geometric representation of the original domain $\Omega$. Compatibility of these displacement fields is requested at the shell mid-surface, such that
\begin{align}
    \vec u_\Bulk|_\Surf &= \vec u_\Surf.
    \label{eq:compatibility}
\end{align}
We note that the restriction $\vec u_\Bulk|_\Surf$ has to be considered in the sense of a trace operation when considering weakly differentiable $\vec u_\Bulk \in [H^1(\Bulk)]^3$. Of course, for the shell thickness $\thickness \to 0$, one expects that $\vec u_\Bulk \to \vec u$. To rigorously prove that 3d elasticity converges to plate and shell models for vanishing thickness is a delicate problem. For linear elasticity Braess et al. proved in \cite{BSS2011} a convergence towards the Reissner--Mindlin and Kirchhoff--Love plates. For (nonlinear) shells results relying on the concept of $\Gamma$-convergence are available, see e.g. \cite{FJGM2003,LMP2010}.

\paragraph{Kinematic reduction}

For thin shells, Kirchhoff--Love kinematic assumptions are often used to describe the through-the-thickness deformation in terms of the shell's mid-surface deformation only \cite{SIMO89,Ciarlet05,BRI17}. Similar assumptions on the displacement field in the shell part $\Omega_\Surf$ are posed in the following. Beforehand, we state
the general representation of a shell surface in the framework of differential geometry, as well as the relevant surface measures. 

The shell's mid-surface $\Surf$ is assumed to be sufficiently smooth to allow for the following definitions of the surface gradient and metric tensors in a classical sense; the extension to surfaces with kinks and T-joints is discussed later in Section~\ref{sec:mixed_shell_formulation}. Let $\nabla_\Surf$ denote the surface gradient, which acts in the tangential plane of the surface only. 
The first and second metric tensor of the shell's mid-surface describe the projection to its tangential plane and its curvature, respectively. 
We use
\begin{align}
    \tsr A &= \tsr I - \vec N_\Surf \vec N_\Surf, &
    \tsr B &= - \nabla_\Surf \vec N_\Surf.
    \label{eq:defAB}
\end{align}
It is noted that, for any smooth field $f$ defined on any three-dimensional superset $\Omega \supset \Surf$, the restriction of its gradient to the surface and projection onto the tangential plane yields the surface gradient, such that
\begin{align}
    \nabla_\Surf f &= \tsr A \cdot (\nabla f)|_{\Surf}. 
\end{align}
Above, the shell's mid-surface displacement in the geometrically reduced problem was introduced as $\vec u_\Surf$. The surface deformation gradient evaluates to,
\begin{align}
    \tsr F_\Surf &= (\nabla_\Surf \vec u_\Surf)^T + \tsr A.
\end{align}
The surface normal to the deformed shell $\vec n_\Surf$ is connected to $\vec N_\Surf$ by \eqref{eq:transformN}, where $\tsr F$ is to be replaced by $\tsr F_\Surf$. 

Kirchhoff--Love kinematic assumptions on $\vec u$ in $\Omega_\Surf$ are now made to connect volumetric and mid-surface displacement; first, it is assumed that the restriction $\vec u|_\Surf = \vec u_\Surf$ holds in the sense of a trace operation. Additionally, any material point $\vec X + \zeta \vec N_\Surf$ is mapped under deformation to the according point on the deformed surface normal,
\begin{align}
    \vec X + \zeta \vec N_\Surf &\to \vec X + \vec u_\Surf(\vec X) + \zeta \vec n_\Surf& & \text{ for }\vec X \in \Surf.
    \label{eq:KL}
\end{align}
The above restriction \eqref{eq:KL} represents the classic conditions of no shear and no thickness deformation.
Note that, in such a regime, the surface deformation gradient $\tsr F_\Surf$ can equivalently be defined as the projection of the volumetric deformation gradient $\tsr F$ to the tangential plane,
\begin{align}
    \tsr F_\Surf &= (\tsr F \cdot \tsr A)|_\Surf.
\end{align}

Finally, surface measures of the deformed surface are given by
\begin{align}
    \tsr a &= \tsr I - \vec n_\Surf \vec n_\Surf, &
    \tsr F_\Surf^T \cdot \tsr b &= - \nabla_\Surf \vec n_\Surf.
    \label{eq:defab}
\end{align}
The Green strain measures, membrane strain $\tsr e$ and curvature $\tsr \kappa$, are defined as the respective changes of the metric tensors after pull-back to the reference configuration. They are symmetric and in the tangential plane, defined through
\begin{align}
    \tsr e &= \frac{1}{2} (\tsr F_\Surf^T \cdot \tsr a \cdot \tsr F_\Surf - \tsr A), &
    \tsr \kappa &= -(\tsr F_\Surf^T \cdot \tsr b \cdot \tsr F_\Surf - \tsr B) = \nabla_\Surf \vec n_\Surf \cdot \tsr F_\Surf - \nabla_\Surf \vec N_\Surf.
    \label{eq:def_ekappa}
\end{align}

Under the Kirchhoff--Love kinematic assumptions \eqref{eq:KL}, the in-plane part of the Cauchy-Green tensor $\tsr C_\Surf = \tsr A \cdot \tsr C \cdot \tsr A$ can be expressed in terms of the strain measures \eqref{eq:def_ekappa} when neglecting higher order terms in the thickness coordinate $\zeta$,
\begin{align}
    \tsr C_\Surf &= \tsr A + 2 \tsr e + 2 \zeta \tsr \kappa +O(\zeta^2).
\end{align}

\paragraph{Reduction of energy densities}
Let us now assume that the free energy of the full problem composes additively into a shell and a bulk part, such that
\begin{align}
    \int_\Omega \psi(\tsr C)\, dV &= \int_{\Omega_\Bulk} \psi_\Bulk(\tsr C)\, dV + \int_{\Omega_\Surf} \psi_\Surf(\tsr C)\, dV.
\end{align}
These respective contributions are to be replaced by approximate energy representations on $\Bulk$ and $\Surf$. For the bulk, it is proposed to use the extended domain $\Bulk$ instead of $\Omega_\Bulk$, amounting to
\begin{align}
    \Psi_\Bulk &= \int_{\Bulk} \psi_\Bulk(\tsr C)\, dV, &
    \delta W_{\text{int},\Bulk} = \delta \Psi_\Bulk = \int_\Bulk \left( \frac{\partial \psi_\Bulk}{\partial \tsr C} : \delta \tsr C \right)\, dV.
    \label{eq:deltaWbulk}
\end{align}

The reduction is less straightforward for the shell domain. 
For certain volume energy densities $\psi_\Surf$, it is possible to derive the in-plane strain-stress relation under plane stress conditions analytically in terms of $\tsr C_\Surf$, and thereby the Green strain measures $\tsr{e}$ and $\tsr{\kappa}$. Let $\tilde \psi_{\Surf}(\tsr C_\Surf)$ denote such a plane-stress continuum energy density for the shell domain. Explicit integration through the thickness then yields a surface energy density $\psi_{\Surf,t}$,
\begin{align}
  \int_{-t/2}^{t/2} \tilde \psi_\Surf(\tsr C_\Surf)\, dV &= \psi_{\Surf,t}(\tsr e, \tsr \kappa).
  \label{eq:defpsis}
\end{align}
In the case of St.~Venant--Kirchhoff materials, one finds
\begin{align}
  \begin{split}
    \psi_{\Surf,t}(\tsr e, \tsr \kappa) &= 
    \frac{t_\Surf\,E}{2(1-\nu^2)}\left( (1-\nu) \tsr e : \tsr e + \nu (\trace\tsr e)^2
    +\frac{t_\Surf^2}{12} \left( (1-\nu) \tsr \kappa : \tsr \kappa + \nu (\trace\tsr \kappa)^2\right)\right)\\
    &= \frac{1}{2} \left( \tsr e : \mathbb C_m : \tsr e + \tsr \kappa : \mathbb C_b : \tsr \kappa\right).
  \end{split}
  \label{eq:defpsisSVK}
\end{align}
Above, $\mathbb C_m$ and $\mathbb C_b$ denote the fourth order tensors representing membrane and bending stiffness for the plane stress case, respectively. It is noted that membrane stiffnesses grow linearly with the thickness parameter, $\mathbb C_m \in O(t_\Surf)$, while bending stiffnesses are of third order, $\mathbb C_b \in O(t_\Surf^3)$.

We propose to use expression \eqref{eq:defpsis} for the shell contribution to the virtual work of internal forces,
\begin{align}
    \Psi_\Surf &= \int_\Surf \psi_{\Surf,t}(\tsr e, \tsr \kappa)\, d\Surf, &
    \delta W_{\text{int},\Surf} = \delta \Psi_\Surf = \int_\Surf \left( \frac{\partial \psi_{\Surf,t}}{\partial \tsr e} : \delta \tsr e + \frac{\partial \psi_{\Surf,t}}{\partial \tsr \kappa} : \delta \tsr \kappa  \right)\, d\Surf.
    \label{eq:deltaWsurf}
\end{align}
It is noted that the energy density $\psi_{\Surf,t}$ was derived under the assumption of a plane-stress state, i.e. for vanishing shear and thickness stresses. This assumption is not necessarily met when structures are immersed in bulk material; for many applications it is valid to assume that shear and thickness stresses are small compared to the in-plane stresses met with in the structural reinforcements. Computational results in Section~\ref{sec:computational_results} support this assumption.

Combining \eqref{eq:deltaWbulk} and \eqref{eq:deltaWsurf} into \eqref{eq:PVA}, we pose the variational problem to find $\vec u_\Bulk$, $\vec u_\Surf$ satisfying the compatibility constraint \eqref{eq:compatibility}, such that for all admissible virtual displacements $\delta \vec u_\Bulk$ and $\delta \vec u_\Surf$ and dependent virtual membrane strains $\delta \tsr e$ and virtual curvatures $\delta \tsr \kappa$, equilibrium is achieved;
\begin{align}
    \int_\Bulk \left( \frac{\partial \psi_\Bulk}{\partial \tsr C} : \delta \tsr C\right)\, dV +
    \int_\Surf \left( \frac{\partial \psi_{\Surf,t}}{\partial \tsr e} : \delta \tsr e + \frac{\partial \psi_{\Surf,t}}{\partial \tsr \kappa} : \delta \tsr \kappa  \right)\, d\Surf &= 
    \int_\Bulk \left( \frac{\partial \psi_\mathrm{ext}}{\partial \vec u} \cdot \delta \vec u \right)\, dV.
    \label{eq:PVA2}
\end{align}
Let us shortly discuss the notion of \emph{admissibility} of displacement fields in the current setting employing classical shell theory. For the second metric tensor $\tsr B$ and curvature tensor $\tsr \kappa$ to be well-defined in the sense of classical derivatives, it is necessary to have the normal vector $\vec N_\Surf$ piecewise $\mathcal{C}^1$ and globally continuous, and $\vec u_\Surf$ to be continuously differentiable and piecewise $\mathcal{C}^2$ \emph{at least}. This is a severe restriction on the underlying geometry -- kinks are not allowed -- as well as on any subsequent finite element discretization. A strategy to alleviate these constraints, and thereby enable the usage of standard nodal elements, is discussed using mixed shell elements.


\paragraph{The case of small deformations and plates}

The equations simplify greatly when restricted to the small deformation case and plates. Assume that not only the usage of the geometrically linear theory of elasticity using the linearized strain $\tsr \eps$ and stress $\tsr \sigma$ according to \eqref{eq:linearstress} can be rectified, but also that $\Surf$ is non-curved. Under these assumptions, we shortly provide the ensuing formulation. For the case of a linear Kirchhoff--Love plate, we assume without restriction of generality that the unit normal is directed in $z$ direction, such that $\vec N_\Surf = \vec e_z$. Consistently, we address the deflection component of $\vec u_\Surf$ by $u_z := \vec u_\Surf \cdot \vec e_z$.
The deformed surface normal is approximated by $\vec n_\Surf = (\tsr I - \nabla_\Surf \vec u_\Surf)\cdot \vec e_z = \vec e_z - \nabla_\Surf u_z$. 
The surface gradient $\nabla_\Surf$ acts in the $\vec e_x \vec e_y$ plane only. Surface metric tensors simplify to $\tsr A = \tsr I - \vec e_z \vec e_z$ the projection into the plane and $\tsr B = \tsr 0$. The linearized surface measures are
\begin{align}
  \tsr e_{\mathrm{lin}} &= \frac{1}{2}(\nabla_\Surf \vec u \cdot \tsr A + \tsr A \cdot \nabla_\Surf \vec u^T), &
  \tsr \kappa_{\mathrm{lin}} &= - \nabla_\Surf^2 u_z.
\end{align}
The variational principle \eqref{eq:PVA2} simplifies greatly for the linearized problem,
\begin{align}
  \int_\Bulk \left( \tsr \eps : \mathbb C : \delta \tsr \eps\right)\, dV +
  \int_\Surf \left( \tsr e_{\mathrm{lin}} : \mathbb C_m : \delta \tsr e_{\mathrm{lin}} + \tsr \kappa_{\mathrm{lin}} : \mathbb C_b : \delta \tsr \kappa_{\mathrm{lin}}  \right)\, d\Surf &= 
  \int_\Bulk \left( \frac{\partial \psi_\mathrm{ext}}{\partial \vec u} \cdot \delta \vec u \right)\, dV.
  \label{eq:PVA2lin}
\end{align}

\section{Consistent discretization using mixed shell elements}
\label{sec:mixed_shell_formulation}
A consistent discretization of \eqref{eq:PVA2} makes it necessary to find displacements $\vec u \in [H^1(\Bulk)]^3$ with the additional constraint of higher smoothness on the shell surface, $\vec u|_\Surf = \vec u_\Surf \in [H^2(\Surf)]^2$. Otherwise, the curvature $\tsr \kappa$ depending on $\vec u_\Surf$ as provided in \eqref{eq:def_ekappa} is not well-defined. This requirement is in direct opposition the usage of low-regularity elements. In the following, we describe the approach first mentioned by Neunteufel and Schöberl \cite{Neunteufel:2019}; a mixed finite element shell formulation which renders this requirement of higher smoothness obsolete. 

Let $\mathcal{T} = \{T\}$ be a compatible finite element mesh of the domain $\Bulk$ that resolves the shell surface $\Surf$ by the two-dimensional (sub-)surface mesh $\mathcal{T}_\Surf = \{T_\Surf\}$. We aim at providing a variant of the principle of virtual work \eqref{eq:PVA2} that is valid when using standard nodal $H^1$ conforming elements for the displacement field $\vec u$, and when setting $\vec u_\Surf = \vec u|_\Surf$ the restriction to the surface mesh $\mathcal{T}_\Surf$. A careful evaluation of derivatives shows that only the contributions due to bending in \eqref{eq:PVA} have to be treated explicitly, all other derivatives are well-defined in weak sense for the proposed continuous finite element displacement functions.
We are therefore concerned with the virtual work of internal forces corresponding to the shell energy \eqref{eq:deltaWsurf}.

\begin{figure}
  \centering
  \begin{tikzpicture}[scale=0.6]

    \draw[line width=1pt, color=gray] (1,2) -- (2,0) -- (3,3);
    \draw[line width=1pt, color=gray] (3,3) -- (5,0.5) -- (7,3);
    \draw (4,2.6) node {$T^-$};
    \draw[line width=1pt, color=gray] (7,3) -- (8,0.5) -- (9.5,2.5);
    \draw[line width=1pt, color=gray]  (2,0) -- (5,0.5) -- (8,0.5);

    \draw[line width=1pt, color=gray] (1,2) -- (1,4) -- (3,3);
    \draw[line width=1pt, color=gray] (3,3) -- (4.5,5.5) -- (7,3);
    \draw (4,3.5) node {$T^+$};
    \draw[line width=1pt, color=gray] (7,3) -- (9,4.5) -- (9.5,2.5);
    \draw[line width=1pt, color=gray]  (1,4) -- (4.5,5.5) -- (9,4.5);

    \draw[line width=1.5pt,-stealth] (5,1.5) -- (5,3);
    \draw (6.5,2.3) node {$\vec N^- = \vec N_\Surf$};

    \draw[line width=1.5pt,-stealth] (5,4.5) -- (5,3);
    \draw (5.6,3.8) node {$\vec N^+$};

    \draw[line width=1.5pt] (1,2) -- (3,3) -- (7,3) -- (9.5,2.5);
    \draw (0.8,1.5) node {$\Surf$};

  \end{tikzpicture}
  \begin{tikzpicture}[scale=0.6]
  
    \draw[line width=1pt] (0,0) -- (4,1);
    \draw[line width=1pt] (0,0) -- (5,4);
    \draw[line width=1pt] (4,1) -- (5,4);
    \draw[line width=1pt] (4,1) -- (10,0);
    \draw[line width=1pt] (5,4) -- (10,0);
    \draw (2,1) node {$T_l$};
    \draw (7.5,1) node {$T_r$};
    \draw (5.1,3.4) node {$E$};

    \draw (4,0.5) node {$V_E^{(1)}$};
    \draw (5,4.5) node {$V_E^{(2)}$};
  
    \draw[line width=1.5pt,-stealth] (6,2) -- (7,4);
    \draw[line width=1.5pt,-stealth] (3,2) -- (2,4);
    \draw (7.2,3.2) node {$\vec N_r$};
    \draw (1.9,3.2) node {$\vec N_l$};
  
    \draw[line width=1.5pt,-stealth] (6,1.5) -- (4.3333,2);
    \draw (5.5,1.3) node {$\vec \nu_r$};
    \draw[line width=1.5pt,-stealth] (3,1.6) -- (4.3333,2);
    \draw (3.2,1.3) node {$\vec \nu_l$};

    \draw[line width=1.5pt,-stealth] (4.3333,2) -- (4.75,3.25);
    \draw (5,2.5) node {$\vec t_E$};
  \end{tikzpicture}
  \caption{Normal and tangential vectors in the finite element mesh. Left: surface normals. Right: in-plane edge normals.}
  \label{fig:normalstangents}
\end{figure}
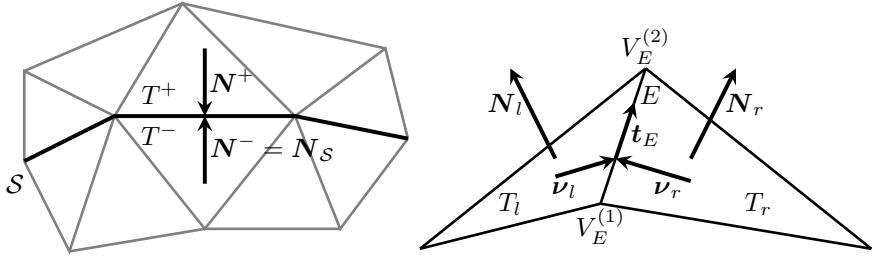

\paragraph{A note on jumps and normal vectors}
On each surface element $T_\Surf \in \mathcal{T}_\Surf$, the surface normal vector in undeformed configuration $\vec N_\Surf$ is uniquely defined up to a sign. In general, such a surface element is shared by two volume elements $T^+$ and $T^-$. Without restriction of generality, we assume that the outward unit normal $\vec N^-$ of $T^-$ coincides with $\vec N_\Surf$, while $\vec N^+ = - \vec N_\Surf$, see also Figure~\ref{fig:normalstangents}. Then, we define the \emph{normal jump} of any tensor- or vector-valued quantity $\tsr \alpha$ across $T_\Surf$ uniquely as
\begin{align}
  [\![\tsr \alpha \cdot \vec N_\Surf]\!]_{T_\Surf} &:= (\tsr \alpha \cdot \vec N_\Surf)|_{T^-} - (\tsr \alpha \cdot \vec N_\Surf)|_{T^+} = (\tsr \alpha^- \cdot \vec N^-) + (\tsr \alpha^+ \cdot \vec N^+).
  \label{eq:jumpF}
\end{align}
An equivalent definition respecting the surface normal in deformed configuration can be made for $[\![\vec \alpha \cdot \vec n_\Surf]\!]_{T_\Surf}$.
If the surface $\Surf$ is part of the boundary $\partial \Omega$, we assume, again without restriction of generality, that the surface normal $\vec N_\Surf$ and boundary outward normal $\vec N$ coincide, such that $T^+$ does not exist. In this case, we formally define the jump as
\begin{align}
  [\![\tsr \alpha \cdot \vec N_\Surf]\!]_{T_\Surf} &:= (\tsr \alpha \cdot \vec N_\Surf)|_{T^-} = \tsr \alpha \cdot \vec N.
  \label{eq:jumpFbd}
\end{align}

Let $\mathcal{E}_\Surf = \{E\}$ denote the set of oriented edges in the surface mesh $\mathcal{T}_\Surf$, also called skeleton of $\mathcal{T}_{\Surf}$.
Each oriented edge $E \in \mathcal{E}_\Surf$ is shared by two triangles $T_l, T_r \in \mathcal{T}_\Surf$. Let now $\vec N_l = \vec N_\Surf|_{T_l}$ and $\vec N_r = \vec N_\Surf|_{T_r}$ denote the surface normals of left and right surface element, respectively. Additionally, $\vec \nu_l$ and $\vec \nu_r$ are the in-plane outward normal vectors of the two triangles, which are only collinear if the discrete shell surface is smooth. The unique edge tangent $\vec t_E$ is oriented pointing from vertex $V_E^{(1)}$ to vertex $V_E^{(2)}$. Without restriction of generality, we assume that the orientation of the edge is such that $T_l$ lies indeed on the left-hand side of edge $E$ as indicated in Figure~\ref{fig:normalstangents}, and set the unique edge normal $\vec \nu = \vec \nu_l$.
For such an edge, we define the \emph{in-plane edge normal jump} of a tensor or vector field $\tsr \tau$ as
\begin{align}
  [\![\tsr \tau \cdot \vec \nu]\!]_{E_\Surf} &:= (\tsr \tau_l \cdot \vec \nu_l) + (\tsr \tau_r \cdot \vec \nu_r).
  \label{eq:jumpE}
\end{align}
Additionally, we need the \emph{jump and average of the in-plane normal-normal component} of a tensor field $\tsr \tau$ on edge $E_\Surf$,
\begin{align}
  [\![ \tau_{\nu\nu} ]\!]_{E_\Surf} &:= [\![ \tsr \tau \cdot \vec \nu]\!]_{E_\Surf} \cdot \vec \nu = (\vec \nu_l \cdot \tsr \tau_l \cdot \vec \nu_l) - (\vec \nu_r \cdot \tsr \tau_r \cdot \vec \nu_r),
  \label{eq:jumpnnE}\\
  \{\!\{ \tau_{\nu\nu} \}\!\}_{E_\Surf} &:= \frac{1}{2}\big((\vec \nu_l \cdot \tsr \tau_l \cdot \vec \nu_l) + (\vec \nu_r \cdot \tsr \tau_r \cdot \vec \nu_r)\big).
  \label{eq:avgnnE}
\end{align}

\paragraph{Low-regularity shell elements}
Neunteufel and Schöberl \cite{Neunteufel:2019} proposed to discretize the moment tensor $\tsr m = \partial \psi_{\Surf,t} / \partial \tsr \kappa$ as a separate unknown field $\hat {\tsr m}$ sporting in-plane normal-normal continuity. To make this approach applicable to general non-linear shell energies $\psi_{\Surf,t}(\tsr e, \tsr \kappa)$, a three-field approximation including another independent field, curvature $\hatkappa$ is convenient, see the recent work \cite{Neunteufel:2023}. In short, let $\mathcal{T}_\Surf = \{T_\Surf\}$ be a surface mesh consisting of (possibly curved) triangular elements $T_\Surf$, and let $\mathcal E_\Surf$ denote the set of element edges. Then, three independent fields $\vec u_\Surf$, $\hatkappa$, and $\hatm$ are discretized through
\begin{subequations}
  \begin{align}
    \vec u_\Surf, \delta \vec u_\Surf &\in \{ \vec v \in [\mathcal{C}(\Surf)]^3: \vec v|_{T_\Surf} \in [P^k(T_\Surf)]^3\}, \label{eq:fe_u0}\\
    \hatkappa, \delta \hatkappa &\in \{\hatkappa \in L^2_{\mathrm{sym}}: \hatkappa|_{T_\Surf} \in P_{\mathrm{sym}}^{k-1}(T_\Surf)\}, \label{eq:fe_kappa}\\
    \hatm, \delta \hatm &\in \{ \hatm \in L^2_{\mathrm{sym}}: \hatm|_{T_\Surf} \in P_{\mathrm{sym}}^{k-1}(T_\Surf), [\![\hat m_{\nu\nu}]\!]_E = 0 \text{ on } E \in \mathcal{E}_\Surf \}. \label{eq:fe_m}
  \end{align}
\end{subequations}
Above, $\vec \nu_r$ and $\vec \nu_l$ are the in-plane normal vectors of the edge $E$ connecting elements $T_r$ and $T_l$ in undeformed configuration, see also Figure~\ref{fig:normalstangents} for a sketch. The set of all symmetric square-integrable matrices in the tangential plane are denoted by $L^2_{\mathrm{sym}}$. For details on the construction of these elements, also on hybridization techniques, we refer to the original work \cite{Neunteufel:2019}.

While the tensor-valued independent moment tensor $\hatm$ acts as a Lagrangian multiplier enforcing the definition of curvature, $\hatkappa = \nabla_\Surf \vec n\cdot \tsr F_\Surf - \nabla_\Surf \vec N$, its characteristic equilibrium property $[\![\hat m_{\nu\nu}]\!]_E = 0$ ensures that $\vec \nu_r \cdot \hatm \cdot \vec \nu_r =\vec  \nu_l \cdot \hatm \cdot \vec \nu_l =: \hat m_{\nu\nu}$ is unique on element edges: on each edge, the \emph{change of angle} between left-hand-side and right-hand-side triangle normals $\vec n_l$ and $\vec n_r$ is dual to the edge moment $\hat m_{\nu\nu}$, and enters into the work statement. The shell energy is enriched by the corresponding terms, defining the potential
\begin{align}
\begin{split}
  \Pi_\Surf(\vec u_\Surf, \hatkappa, \hatm) := &
  \int_{\Surf} \psi_{\Surf,t}(\tsr e, \hatkappa)\, d\Surf
  + \sum_{T_\Surf \in \mathcal{T}_\Surf} \int_{T_\Surf} (-\hatkappa + \nabla_\Surf \vec n \cdot \tsr F_\Surf - \nabla_\Surf \vec N_\Surf) : \hatm\, d\Surf\, +\\
  & \sum_{E \in \mathcal E_\Surf}\int_E (\operatorname{acos}(\vec n_l \cdot \vec n_r)-\operatorname{acos}(\vec N_l \cdot \vec N_r)) \hat m_{\nu\nu}\, ds \to \min_{\vec u_\Surf,\hatkappa \text{ adm.}}\max_{\hatm \text{ adm.}}\,.
\end{split}
\label{eq:defPi}
\end{align}

The original shell energy $\Psi_\Surf$ from \eqref{eq:deltaWsurf} is replaced by $\Pi_\Surf$ in the principle of virtual work. Finally, for all kinematically admissible virtual displacements $\delta \vec u$, virtual curvatures $\delta \hatkappa$ and virtual moment tensors $\delta \hatm$ as in \eqref{eq:fe_m} and \eqref{eq:fe_m}, the variational equation is satisfied,
\begin{align}
 \delta \Psi_\Bulk + \delta \Pi_\Surf &= \int_\Bulk \left( \frac{\partial \psi_\mathrm{ext}}{\partial \vec u} \cdot \delta \vec u \right)\, dV.
\label{eq:variational_nonlin}
\end{align}
Note that the explicit computation of the variation $\delta \Pi_\Surf$ includes the variation $\delta \vec n$. Computing $\delta \vec n$ for the current displacement state $\vec u$ and a given virtual displacement $\delta \vec u$ is a tedious task if done by hand. In implementations, automatic differentiation can treat these nonlinear dependencies easily. As demonstrated by Neunteufel and Schöberl in \cite{Neunteufel:2021} the shell element is free of membrane locking when interpolating the membrane strains $\tsr e$ locally into the so-called Regge space. We refer to Christiansen \cite{Chr2011} and Li \cite{li18} for the development of Regge finite elements from Regge calculus \cite{Regge61}, but we do not go into the details here.

\paragraph{Hybridization and structures with kinks and branches}
The shell formulation \eqref{eq:defPi} has a saddle-point structure, which leads to an indefinite stiffness matrix after assembling. To overcome this issue, hybridization techniques can be used, see e.g. \cite{Hughes00,BBF13} and specifically for the TDNNS method \cite{Neunteufel:2023,Neunteufel:2021threefield} for details. The idea is to break the normal-normal continuity of the moment tensor elements \eqref{eq:fe_m}. By adding an additional unknown as Lagrange multiplier $\vec{\alpha}$, we reinforce the broken continuity in weak sense. Therefore, the moment tensor $\hatm$ (and the curvature tensor $\hatkappa$) can be efficiently eliminated from the system matrix element-wise, and the resulting system involving only the displacement fields $\vec u$ and the new Lagrange multiplier $\vec\alpha$ is a minimization problem.

The additional unknowns live only on the skeleton $\mathcal{E}_\Surf$, the edges, of the shell triangulation $\mathcal T_\Surf$ and will be denoted by $\vec{\alpha}$.
The set of all admissible Lagrange functions is given by
\begin{align*}
  \vec{\alpha} \in \{ \vec{\beta} \in [\mathcal{C}(\mathcal{E}_\Surf)]^3: \vec{\beta}|_{E} \in [P^k(E)]^3,\, [\![\vec \beta \cdot \vec \nu]\!]_{E} = \vec 0,\, \vec N_\Surf\cdot \vec{\beta}=\vec t_E\cdot \vec{\beta}=0, \vec{\beta}=\vec{0} \text{ on }\Gamma_{\text{fix}}\}.
\end{align*}
This choice implies that for any $\vec \alpha$ as above, the in-plane normal component $\alpha_\nu := \vec \alpha_l \cdot \vec \nu_l = -\vec \alpha_r \cdot \vec \nu_r$ is unique. The hybridized version of \eqref{eq:defPi} reads
\begin{align}
\begin{split}
  \Pi_\Surf(&\vec u_\Surf, \hatkappa, \hatm, \vec{\alpha}) := 
  \int_{\Surf} \psi_{\Surf,t}(\tsr e, \hatkappa)\, d\Surf
  + \sum_{T_\Surf \in \mathcal{T}_\Surf} \int_{T_\Surf} (-\hatkappa + \nabla_\Surf \vec n \cdot \tsr F_\Surf - \nabla_\Surf \vec N_\Surf) : \hatm\, d\Surf\, +\\
  & \sum_{E \in \mathcal E_\Surf}\int_E ( \operatorname{acos}(\vec n_l \cdot \vec n_r)-\operatorname{acos}(\vec N_l \cdot \vec N_r))\, \{\!\{\hat m_{\nu\nu}\}\!\}_E+\vec{\alpha}_{\nu}[\![\hat m_{\nu\nu}]\!]\, ds\to \min_{\vec u_\Surf, \hatkappa,\vec\alpha \text{ adm.}}\max_{\hatm \text{ adm.}}\,,
\end{split}
\label{eq:defPiHybrid}
\end{align}
where jump and average of the in-plane normal-normal component $\hat m_{\nu\nu}$ are defined in \eqref{eq:jumpnnE}--\eqref{eq:avgnnE}.

For structures with kinks the shell method intrinsically preserves the angle of the kink and the inflow of the normal-normal component of the moment tensor $\hat{m}_{\nu\nu}$ equals the outflow. The moment conservation follows by the construction of $\hatm$. To see that the angle is preserved in a weak sense, consider on edge $E\in\mathcal{E}_{\Surf}$. Then the variation of $\hat m_{\nu\nu}$ at the edge terms of \eqref{eq:defPi} reads 
\begin{align}
  \label{eq:weakangle}
  \int_E (\operatorname{acos}(\vec n_l \cdot \vec n_r)-\operatorname{acos}(\vec N_l \cdot \vec N_r))\, \delta\hat m_{\nu\nu}\, ds=0
\end{align}
and thus $(\operatorname{acos}(\vec N_l \cdot \vec N_r) - \operatorname{acos}(\vec n_l \cdot \vec n_r))=0$.

Branched shells, where more than two elements $\{T_i\}_{i=1}^k$ share a common edge $E\in\mathcal E_\Surf$, are more delicate \cite{Neunteufel:2023}. First, there is no notation of a ``left'' or ``right'' element, such that we have to rewrite the boundary terms into an equivalent form applicable for those structures. To this end, we define the averaged normal vector $\{\!\{\vec N\}\!\}=\frac{\vec N_l+\vec N_r}{|\vec N_l+\vec N_r|}$ and analogously $\{\!\{\vec n\}\!\}$. Then, one can readily show that
\begin{align*}
  \sum_{E \in \mathcal E_\Surf} (\operatorname{acos}(\vec n_l \cdot \vec n_r)-\operatorname{acos}(\vec N_l \cdot \vec N_r))\,\hat m_{\nu\nu} = \sum_{T\in\mathcal{T}_\Surf} (\operatorname{acos}(\vec n_{T}\cdot \{\!\{\vec n\}\!\})-\operatorname{acos}(\vec N_{T}\cdot \{\!\{\vec N\}\!\}))\,\hat m_{\nu\nu}.
\end{align*}
This concept, however, can now be easily extended to branched shells by defining the averaged normal vector $\{\!\{\vec N\}\!\}=\frac{\sum_{i=1}^k\vec N_{T_i}}{|\sum_{i=1}^k\vec N_{T_i}|}$ (to avoid cancellation of the normals one can swap the orientation of one branch if necessary). One important difference to kinked shells is that we require to use the hybridized version, the other method fails to handle branched shells, which can be seen as follows: If we keep $\hatm$ to be normal-normal continuous we would analogously to \eqref{eq:weakangle} preserve the \emph{sum} of all angles. Breaking the continuity of $\hat{m}_{\nu\nu}$ yields the constraints that the angle $\operatorname{acos}(\vec N_{T_i}\cdot \{\!\{\vec N\}\!\})$ is preserved in a weak sense for all $i=1,\dots,k$, i.e., each individual angle is conserved. The hybridization variable $\vec\alpha$ is still single valued on $E$. Instead of the normal-normal continuity of $\hatm$ at the branched edge, it now enforces the \emph{balance of moments}: the sum of all inflow and outflow moments at the edge are equal.  

\paragraph{Linear plate elements}

The variational principle \eqref{eq:variational_nonlin} simplifies greatly for the linearized problem on plates. There holds
\begin{align}
  \delta \Psi_\Bulk + \delta \Pi_{\Surf,\mathrm{lin}} &= \int_\Omega \left( \frac{\partial \psi_\mathrm{ext}}{\partial \vec u} \cdot \delta \vec u \right)\, dV,
 \label{eq:variational_lin}
 \end{align}
with
\begin{align}
  \begin{split}
  \Pi_{\Surf,\mathrm{lin}} \ = &
  \int_\Surf \left( \frac{1}{2}\tsr e_{\mathrm{lin}} : \mathbb C_m : \tsr e_{\mathrm{lin}} + \frac{1}{2} \tsr \kappa_{\mathrm{lin}} : \mathbb C_b : \tsr \kappa_{\mathrm{lin}}  - \tsr m_{\mathrm{lin}} : \tsr \kappa_{\mathrm{lin}} \right)\, d\Surf - \\
   & \sum_{T_\Surf \in \mathcal{T}_\Surf}  
    \int_{T_\Surf} \nabla_\Surf^2 u_z : \tsr m_{\mathrm{lin}} \, d\Surf\, - \sum_{E \in \mathcal{E}_\Surf}
  \int_{E} [\![\nabla_\Surf u_z \cdot \vec \nu]\!]   m_{\mathrm{lin},\nu\nu} \, ds. 
  \end{split}
  \label{eq:Pi_S_lin}
\end{align}
The derivation of the above terms can be found in similar form in \cite{Neunteufel:2019}. Of course, hybridization can be included in an equivalent manner. When eliminating $\tsr \kappa_{\mathrm{lin}}$ from \eqref{eq:Pi_S_lin} the mixed two-field Hellan-Herrmann-Johnson method \cite{Hel67,Her67,Joh73} for Kirchhoff--Love plates is recovered.

\section{Implications on local equilibrium}
\label{sec:local_equilibrium}
In this section, we restrict the analysis to the case of linear elasticity, small deformation, and non-curved plates.

Under the above assumptions on the energy approximations and kinematic restrictions for admissible displacement fields, the principle of virtual work can be used to deduce local equilibrium equations. Within this section, we compute the local equilibrium conditions at the shell/plate mid-surface, and discuss their validity.

In classical continuum mechanics, local equilibrium of forces at small deformations is characterized in differential form, $\opdiv \tsr \sigma + \vec f = 0$ for volume forces $\vec f$, and on the boundary $\tsr \sigma \cdot \vec N = \vec t$ for the surface traction forces $\vec t$. The stress tensor is not necessarily differentiable everywhere in classical sense; however, on any interface within the body, with $\vec N$ its normal vector, continuity of the stress vector $\tsr \sigma \cdot \vec N$ goes back to Cauchy. In any classical displacement-based finite element discretization, continuity of the stress vector $\tsr \sigma \cdot \vec N$ is enforced in weak sense on element interfaces, such that any discontinuity does not perform work on the discrete virtual displacements.
This is usually derived from the principle of virtual work, performing integration by parts on each element separately. As similar considerations will be needed for the much more complex case of the coupled continuum/shell model, we will perform this analysis in detail for the well-understood case in classical continuum theory.

\paragraph{Local equilibrium in linear elasticity problems}
Let now $\mathcal{T} = \{T\}$ be a compatible finite element mesh of a domain $\Omega$. 
 Element-wise integration by parts in the virtual work of internal forces leads to
\begin{align}
  \int_\Omega \tsr \sigma : \delta \tsr \eps \, dV &=
  \sum_{T \in \mathcal T} \left(\int_T -\opdiv \tsr \sigma \cdot \delta \vec u\, dV + \int_{\partial T} (\tsr \sigma \cdot \vec N_T) \cdot \delta \vec u\, dA\right).
\end{align}
Reordering the surface integrals face by face, using definition \eqref{eq:jumpF} of the jump of $\tsr \sigma \cdot \vec N$ across $F \in \mathcal F$, where $\mathcal F$ denotes the set of all facets of $\mathcal T$, this is equivalent to
\begin{align}
  \int_\Omega \tsr \sigma : \delta \tsr \eps \, dV &=
  \sum_{T \in \mathcal T} \int_T-\opdiv \tsr \sigma \cdot \delta \vec u\, dV + \sum_{F \in \mathcal F} \int_{F} [\![\tsr \sigma \cdot \vec N]\!]_F \cdot \delta \vec u\, dA.
\end{align}
For simplicity, assume that only volume forces $\vec f$, no surface tractions $\vec t$ are acting on the domain. Then, the principle of virtual work states
\begin{align}
  \sum_{T \in \mathcal T} \int_{T} (-\opdiv \tsr \sigma + \vec f) \cdot \delta \vec u\, dV + \sum_{F \in \mathcal F} \int_{F} [\![\tsr \sigma \cdot \vec N]\!]_F \cdot \delta \vec u\, dA = 0.
\end{align}
In other words, equilibrium of forces $-\opdiv \tsr \sigma + \vec f$ and continuity of the stress vector $[\![\tsr \sigma \cdot \vec N]\!]_F = \vec 0$ are satisfied in weak sense.

\paragraph{Local equilibrium in coupled plate problems}
In the following, a similar derivation will be done for the coupled linearized problem as stated in \eqref{eq:variational_lin}. To simplify the derivation, we first assume that 
the displacement field's derivative is globally continuous. Then the edge terms in the linear plate formulation \eqref{eq:Pi_S_lin} vanish. Further, we assume that $\tsr m_{\mathrm{lin}}$ is continuous at vertices and that no external forces are acting, such that the work of external forces equals zero, $\delta W_{\mathrm{ext}} = 0$. Otherwise, external forces enter the equilibrium conditions in a well-known manner.

Define, for ease of presentation $\tsr{\tau}:=\mathbb{C}_m\tsr{e}_{\mathrm{lin}}$ as the membrane stress tensor and use the symmetric stress tensor $\tsr \sigma$ for $\tsr P$, as we are in the linear regime. Considering some general virtual displacement field $\delta \vec u$ while setting $\delta \tsr \kappa_{\mathrm{lin}} = \tsr 0$ and $\delta \tsr m_{\mathrm{lin}} = \tsr 0$, the principle of virtual work \eqref{eq:variational_lin} reduces to
\begin{align}
 \begin{split}
   &\int_{\Bulk} \tsr \sigma : \delta \tsr \eps\, dV +
    \int_{\Surf} \tsr{\tau}:\delta\tsr{e}_{\mathrm{lin}}\,d\Surf -
     \int_{\Surf} \nabla_\Surf^2 \delta u_z : \tsr m_{\mathrm{lin}} \, d\Surf
   = 0.
 \end{split}
\end{align} 

The variation of the bulk energy has already been discussed above. We continue with the variation of the membrane part. Here, integration by parts is conducted on each surface element, where again the in-plane edge normal vector $\vec \nu$ enters the boundary expression. A reordering of edge terms leads to the following formula, making use of the jump of the membrane stress vector as defined in \eqref{eq:jumpE},
\begin{align}
\begin{split}
     \int_\Surf \tsr \tau : \delta \tsr e_{\mathrm{lin}}\, d\Surf
    &= \sum_{T_\Surf \in \mathcal T_\Surf} \left(-\int_{T_\Surf} \opdiv_\Surf \tsr \tau \cdot \delta \vec u\, d\Surf + \int_{\partial T_\Surf} \tsr \tau \cdot \vec \nu \cdot \delta \vec u\, ds\right) \\
    &=
    \sum_{T_\Surf \in \mathcal T_\Surf} \int_{T_\Surf} (-\opdiv_\Surf \tsr \tau \cdot \delta \vec u)\, d\Surf + \sum_{E \in \mathcal{E}_\Surf} \int_{E_\Surf} [\![\tsr \tau \cdot \vec \nu]\!]_{E_\Surf} \cdot \delta \vec u\, ds.
\end{split}
\end{align}
Note that, as $\tsr \tau$ as well as $\vec \nu$ lie in the tangential plane, only tangential mid-surface displacement components are affected.

Next, the variation of the bending terms are considered. We start by performing integration by parts, using again that $\tsr m_{\mathrm{lin}}$ acts in the tangential plane only, and that $\tsr m_{\mathrm{lin}}$ is normal-normal continuous
\begin{align}
\begin{split}
  \int_{\Surf} \nabla_\Surf^2 \delta u_z : \tsr m_{\mathrm{lin}} \, d\Surf\, &= \sum_{T\in\mathcal{T}_\Surf}\left(- \int_{T_\Surf} \nabla_\Surf \delta u_z \cdot \opdiv_\Surf \tsr m_{\mathrm{lin}} \, d\Surf\, + \int_{\partial T_\Surf} \nabla u_z\cdot (\tsr m_{\mathrm{lin}})_{\nu}\,ds\right)\\
  & = \sum_{T\in\mathcal{T}_\Surf} \int_{T_\Surf} (-\nabla_\Surf \delta u_z \cdot \opdiv_\Surf \tsr m_{\mathrm{lin}}) \, d\Surf\, + \sum_{E_{\Surf}\in\mathcal{E}_\Surf}\int_{E_{\Surf}} \frac{\partial u_z}{\partial \vec t}[\![ (m_{\mathrm{lin}})_{\nu t}]\!]_{E_{\Surf}}\,ds.
\end{split}
\end{align}
Above, $\partial/\partial \vec t$ denotes the tangential derivative along edges, and $(m_{\mathrm{lin}})_{\nu t}$ is the tangential component of the moment vector $\tsr m_{\mathrm{lin}} \cdot \vec \nu$.
Another integration by parts yields that the virtual work contribution of the bending moment  transforms to 
\begin{align}
\begin{split}
    & \int_{\Surf} \nabla_\Surf^2 \delta u_z : \tsr m_{\mathrm{lin}} \, d\Surf\,  = \sum_{T_\Surf\in\mathcal{T}_\Surf}\int_{T_\Surf} \delta u_z \, \opdiv_\Surf \opdiv_\Surf \tsr m_{\mathrm{lin}} \, d\Surf\,-  \\
    &\sum_{E_{\Surf}\in\mathcal{E}_\Surf}\int_{E_{\Surf}} \delta u_z\left[\!\left[ \frac{\partial (m_{\mathrm{lin}})_{\nu t}}{\partial \vec t}+\vec \nu \cdot \opdiv_\Surf \tsr m_{\mathrm{lin}}\right]\!\right]_{E_{\Surf}}\,ds.
\end{split}
 \end{align}
 Finally, all contributions are collected within a single statement:
\begin{align}
  \begin{split}
    &\sum_{T \in \mathcal T} \int_{T} -\opdiv \tsr \sigma \cdot \delta \vec u\, dV 
    + \sum_{F \in \mathcal F\backslash \mathcal{T}_\Surf} \int_{F} [\![\tsr \sigma \cdot \vec N]\!]_F \cdot \delta \vec u\, dA\, + \\
    &\sum_{T_\Surf \in \mathcal{T}_\Surf} \int_{T_\Surf} \delta \vec u \cdot \left( [\![\tsr \sigma \cdot \vec N]\!]_{T_\Surf} - \opdiv_\Surf \tsr \tau - (\opdiv_\Surf \opdiv_\Surf \tsr m_{\mathrm{lin}})\vec N \right)\, d\Surf\, +\\
    &\sum_{E_\Surf \in \mathcal{E}_\Surf} \int_{E_\Surf} \delta \vec u \cdot \left[\!\left[\tsr \tau \cdot \vec \nu +(\vec \nu \cdot \opdiv_\Surf \tsr m_{\mathrm{lin}} + \frac{\partial  (m_{\mathrm{lin}})_{\nu t}}{\partial \vec t})\,\vec e_z \right]\!\right]_{E_\Surf} \, ds  = 0.
  \end{split}
\end{align} 
Relaxing again the $C^1$ assumptions on the mesh, displacement field, and additional assumptions on the moment tensor one obtains (see \cite{Braess:2020} for a detailed derivation, though in a different context, and \cite[Lemma~5.2]{Braess:2020} for an exact definition of the vertex jump $[\![\cdot]\!]_{V}$)
\begin{align}
  \begin{split}
    &\sum_{T \in \mathcal T} \int_{T} -\opdiv \tsr \sigma \cdot \delta \vec u\, dV 
    + \sum_{F \in \mathcal F\backslash \mathcal{T}_\Surf} \int_{F} [\![\tsr \sigma \cdot \vec N]\!]_F \cdot \delta \vec u\, dA\, + \\
    &\sum_{T_\Surf \in \mathcal{T}_\Surf} \int_{T_\Surf} \delta \vec u \cdot \left( [\![\tsr \sigma \cdot \vec N]\!]_{T_\Surf} - \opdiv_\Surf \tsr \tau - (\opdiv_\Surf \opdiv_\Surf \tsr m_{\mathrm{lin}})\vec N \right)\, d\Surf\, +\\
    &\sum_{E_\Surf \in \mathcal{E}_\Surf} \int_{E_\Surf} \delta \vec u \cdot \left[\!\left[\tsr \tau \cdot \vec \nu +(\vec \nu \cdot \opdiv_\Surf \tsr m_{\mathrm{lin}} + \frac{\partial  (m_{\mathrm{lin}})_{\nu t}}{\partial \vec t})\,\vec e_z \right]\!\right]_{E_\Surf} \, ds\,  - \\
    &\sum_{V_\Surf \in \mathcal{V}_\Surf } \delta u_z(V_\Surf) [\![ (m_{\mathrm{lin}})_{\nu t}]\!]_{V_\Surf}
    = 0.
  \end{split}
\end{align} 

By splitting the test function into normal and tangential parts for the plate terms, i.e., by first setting $\delta\vec u = \delta u_z\,\vec e_z$ and then $\delta \vec u = \delta u_x \vec e_x + \delta u_y \vec e_y$, we obtain the following conditions at vertices $V_\Surf$ and edges $E_\Surf$ of the plate triangulation  $\mathcal{T}_\Surf$, triangles $T_\Surf$ of $\mathcal{T}_\Surf$, all other faces $F$ of the bulk domain, and on volume elements $T$ in the bulk domain, which are discussed in detail below,
\begin{subequations}
  \begin{align}
    &[\![ (m_{\mathrm{lin}})_{\nu t}]\!]_{V_\Surf} = 0,&&\qquad  \forall V_{\Surf} \in \mathcal{V}_{\Surf},\label{eq:cond_vert}\\
    & \left[\!\left[\tsr \tau \cdot \vec \nu  \right]\!\right]_{E_\Surf}=0,\qquad \left[\!\left[\vec \nu \cdot \opdiv_\Surf \tsr m_{\mathrm{lin}} + \frac{\partial (m_{\mathrm{lin}})_{\nu t}}{\partial \vec t} \right]\!\right]_{E_\Surf}=0, &&\qquad  \forall E_{\Surf}\in\mathcal{E}_\Surf,\label{eq:cond_edge}\\
    &\opdiv_\Surf \tsr \tau = [\![\tsr \sigma_{N \Surf}]\!]_{T_{\Surf}}, \qquad \opdiv_\Surf \opdiv_\Surf \tsr m_{\mathrm{lin}} = [\![\sigma_{NN}]\!]_{T_{\Surf}}, &&\qquad \forall T_{\Surf}\in\mathcal{T}_{\Surf},\label{eq:plate_trig}\\
    & [\![\tsr \sigma \cdot \vec N]\!]_{F} = 0,&& \qquad \forall F \in \mathcal F\backslash \mathcal{T}_\Surf,\label{eq:cond_bulk_face}\\
    &-\opdiv \tsr \sigma=0,&& \qquad \forall T \in \mathcal T. \label{eq:cond_tet}
  \end{align}
\end{subequations}

Condition \eqref{eq:cond_vert} enforces a jump condition on the vertices for the moment stress tensor, stating vanishing residual forces, known from the theory of Kirchhoff--Love plates \cite{VNS2007}. On the edges of the plate \eqref{eq:cond_edge} gives continuity of the membrane stress and effective transverse shear force in weak sense. The coupling of stresses between the plate and the bulk are given by conditions \eqref{eq:plate_trig}. Here, the stress vector is split into a normal and tangential component, $\tsr \sigma \cdot \vec N = \vec \sigma_{N\Surf} + \sigma_{NN} \vec N_\Surf$. The membrane forces $\opdiv_\Surf \tsr \tau$ have to be in balance with the change in shear stress $[\![\tsr \sigma_{N \Surf}]\!]_{T_{\Surf}}$. The second identity states that the change in normal forces when crossing the shell surface is balanced by moment forces $\opdiv_\Surf \opdiv_\Surf \tsr m_{\mathrm{lin}}$. Conditions \eqref{eq:cond_bulk_face} and \eqref{eq:cond_tet} are the same as for a pure elasticity problem.

Setting $\delta\vec u= \vec 0$ and $\delta \tsr m_{\mathrm{lin}} = \tsr 0$ we directly deduce from the principle of virtual work \eqref{eq:variational_lin} the stress strain relation 
\begin{align}
  \mathbb C_b \tsr \kappa_{\mathrm{lin}} = \tsr m_{\mathrm{lin}}.
\end{align}

Setting $\delta\vec u= \vec 0$ and $\delta \tsr \kappa_{\mathrm{lin}} = \tsr 0$ and first assuming that the displacement field's derivative on it is globally continuous, \eqref{eq:variational_lin} yields the definition of the linearized curvature tensor
\begin{align}
  \label{eq:curv_lin}
  \tsr \kappa_{\mathrm{lin}} = -\nabla_\Surf^2 u_z
\end{align}
as the moment tensor acts as Lagrange multiplier for \eqref{eq:curv_lin}. For a general displacement field $\vec u$ additionally the normal jump of its derivative over edges has to be considered, leading to the relation
\begin{align}
  \int_{\Surf}\tsr \kappa_{\mathrm{lin}}:\delta\tsr m_{\mathrm{lin}}\, d\Surf = \sum_{T\in\mathcal T}\int_T-\nabla_\Surf^2 u_z : \delta \tsr m_{\mathrm{lin}}\, d\Surf - \sum_{E\in\mathcal E_\Surf}\int_E [\![\nabla_\Surf u_z \cdot \vec \nu]\!] \delta m_{\mathrm{lin},\nu\nu}\, ds.
\end{align}
In case of a linear displacement field $\vec u\in [P^1]^3$, the surface Hessian $\nabla_\Surf^2 u_z$ is zero on each triangle such that only the edge terms account for defining the curvature $\tsr \kappa_{\mathrm{lin}}$ similar to the Hellan--Herrmann--Johnson method \cite{Com89,Hel67,Her67,Joh73}.

\paragraph{Extension to shells}
The above presented derivation can be performed also for linear curved shells, see \cite{Neunteufel:2023} for the formulation. Due to the non-constant normal vector at the initial configuration, however, additional terms involving derivatives of the normal vector arise, and are therefore not performed here. For nonlinear shells the situation is even more complex and out of scope of this work.


\section{Computational results}
\label{sec:computational_results}
All numerical examples are performed in the open source finite element library NGSolve\footnote{\href{www.ngsolve.org}{www.ngsolve.org}} \cite{Sch97,Sch14}.

\subsection{Square laminate block}
\label{subsec:square_laminate_block} 

\begin{figure}[ht]
  \centering
  \begin{tikzpicture}[scale=0.5]
    \def\L{10} 
    \def\H{2} 
    \def\tS{0.05} 
    \def\offset{\H/4} 

    \draw[fill=gray!80] (-\L/2,-\H/2,\L/2) -- (0,-\H/2,\L/2) -- (0,\H/2,\L/2) -- (-\L/2,\H/2,\L/2) -- cycle;
    \draw[fill=gray!80] (0,-\H/2,\L/2) -- (\L/2,-\H/2,\L/2) -- (\L/2,\H/2,\L/2) -- (0,\H/2,\L/2) -- cycle;
    \draw[fill=gray!50] (-\L/2,\H/2,-\L/2) -- (0,\H/2,-\L/2) -- (0,\H/2,\L/2) -- (-\L/2,\H/2,\L/2) -- cycle;
    \draw[fill=gray!20] (0,\H/2,-\L/2) -- (\L/2,\H/2,-\L/2) -- (\L/2,\H/2,\L/2) -- (0,\H/2,\L/2) -- cycle;
    \draw[fill=gray!20] (\L/2,\H/2,-\L/2) -- (\L/2,-\H/2,-\L/2) -- (\L/2,-\H/2,\L/2) -- (\L/2,\H/2,\L/2) -- cycle;
    \draw[dashed] (-\L/2,-\H/2,0) -- (\L/2,-\H/2,0) -- (\L/2,\H/2,0) -- (-\L/2,\H/2,0) -- cycle;
    \draw[dashed] (-\L/2,-\H/2,\L/2) -- (-\L/2,-\H/2,-\L/2) -- (-\L/2,\H/2,-\L/2);
    \draw[dashed] (0,-\H/2,\L/2) -- (0,-\H/2,-\L/2) -- (0,\H/2,-\L/2);
    \draw[dashed] (-\L/2,-\H/2,-\L/2) -- (\L/2,-\H/2,-\L/2);

    \draw[-stealth] (0,0,0) -- (0,0,0.9*\L) node[anchor=east]{$x$};   
    \draw[-stealth] (0,0,0) -- (0.8*\L,0,0) node[anchor=south]{$y$}; 
    \draw[-stealth] (0,0,0) -- (0,2*\H,0) node[anchor=west]{$z$};

    \draw[] (\L/2,0,0) node[anchor=south west] {$A$};
    \node[] at (\L/2,0,0){$\bullet$};
    \node[] at (\L/2,\H/2,0){$\bullet$};
    \draw[] (\L/2,\H/2,0) node[anchor=south] {$B$};

    \node[anchor=north] at (1.2*\L/2,-1.5*\H,0) {$\vec u$ prescribed};
    \draw (\L/2,0,\L/2) -- (1.2*\L/2,-1.5*\H,0) -- (\L/2,0,-\L/2);

    \draw[<->] (-\L/2,-1.*\H,\L/2) -- (\L/2,-1.*\H,\L/2) node[midway, below] {$L=100$ mm};
    \draw[<->] (-\L/2*1.1,-\H/2,\L/2) -- (-\L/2*1.1,\H/2,\L/2) node[midway, left] {$H=20$ mm};
\end{tikzpicture}
\begin{tikzpicture}[scale=0.5]
  \draw[line width=1.2,fill=gray!20] (1,2) -- (3,2) -- (3,-2) -- (1,-2);
  \draw[line width=1.2pt,fill=gray!20] (1,2.2) -- (3,2.2) -- (3,1.8) -- (1,1.8);
  \draw[line width=1.2pt,fill=gray!20] (1,-2.2) -- (3,-2.2) -- (3,-1.8) -- (1,-1.8);
  \draw[line width=1,dashed] (1,2) -- (3,2) -- (3,-2) -- (1,-2);
  \draw[<->] (3.5,-2) -- (3.5,0) node[midway, right] {$H/2$};
  \draw[<->] (3.5,0) -- (3.5,2) node[midway, right] {$H/2$};
  \draw[<->] (0.5,1.8) -- (0.5,2.2) node[midway, left] {$t_\Surf$};
  \draw[<->] (0.5,-1.8) -- (0.5,-2.2) node[midway, left] {$t_\Surf$};
  \draw[-stealth] (0,0) -- (6,0) node[anchor=south west] {$x$};

  \def\offs{-5};
  \draw[line width=1.2,fill=gray!20] (1,2+\offs) -- (3,2+\offs) -- (3,-2+\offs) -- (1,-2+\offs);
  \draw[line width=1.2pt,fill=gray!30] (1,1.2+\offs) -- (3,1.2+\offs) -- (3,0.8+\offs) -- (1,0.8+\offs);
  \draw[line width=1.2pt,fill=gray!30] (1,-1.2+\offs) -- (3,-1.2+\offs) -- (3,-0.8+\offs) -- (1,-0.8+\offs);
  \draw[line width=1,dashed] (1,1+\offs) -- (3,1+\offs) -- (3,-1+\offs) -- (1,-1+\offs);
  \draw[<->] (3.5,-1+\offs) -- (3.5,0+\offs) node[midway, right] {$H/4$};
  \draw[<->] (3.5,-0+\offs) -- (3.5,1+\offs) node[midway, right] {$H/4$};
  \draw[<->] (3.5,-2+\offs) -- (3.5,-1+\offs) node[midway, right] {$H/4$};
  \draw[<->] (3.5,1+\offs) -- (3.5,2+\offs) node[midway, right] {$H/4$};
  \draw[<->] (0.5,0.8+\offs) -- (0.5,1.2+\offs) node[midway, left] {$t_\Surf$};
  \draw[<->] (0.5,-0.8+\offs) -- (0.5,-1.2+\offs) node[midway, left] {$t_\Surf$};
  \draw[-stealth] (0,0+\offs) -- (6,0+\offs) node[anchor=south west] {$x$};

\end{tikzpicture}

  \caption{Square laminate plate: geometric setup of the symmetric block with kinematically constrained surfaces indicated in dark gray (left); different variants of through-the-thickness position of reinforcements (right).}
  \label{fig:ex1_geometry}
\end{figure}
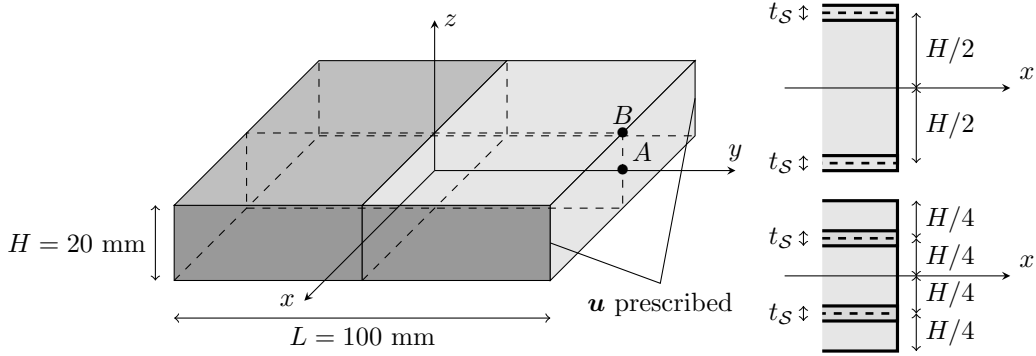

In the first example, a square laminate block is considered. The lateral dimension of the block is chosen as $L = 100$~mm, its bulk height as $H = 20$~mm. The block is reinforced by two layers of high stiffness and thickness $t_\Surf = 0.5$~mm; in a first setup, these layers are located at the bulk's top and bottom surface, while in a second setup the two layers are immersed within the bulk material at a distance of $\pm H/4$ to the block's mid-surface. The block is assumed to be aligned with the spatial coordinate system as depicted in Figure~\ref{fig:ex1_geometry}. 

Both the reinforcing plate and the bulk material are described through the St.~Venant--Kirchhoff model, with Young's modulus $E_\Surf = 10$~GPa and Poisson ratio $\nu_\Surf = 0.3$ for the reinforcement part, and varying Young's modulus $E_\Bulk$ and $\nu_\Bulk = 0.3$ in the bulk material. 
Two test cases prescribing extension in $x$ direction as well as simple shear deformation at two opposite lateral block faces are monitored for different ratios $E_\Bulk/E_\Surf \in [0.001, 1]$. 
In both cases, the block is kinematically free on its top and bottom face, as well as two of the lateral faces; $\tsr P \cdot \vec e_z = \vec 0$ on the top and bottom surface as well as $\tsr P \cdot \vec e_y = \vec 0$ on the lateral surfaces $y = \pm L/2$. Concerning the surfaces at $x = \pm L/2$, kinematic boundary conditions for $\vec u$ are prescribed; for the first load case referred to as \emph{extension}, we set $\vec u = \pm \Delta L/2 \, \vec e_x$ resembling a prescribed normal strain of $\eps_{xx} = \Delta L/L$, for the second load case referred to as \emph{simple shear}, setting $\vec u = \pm \Delta L/2 \vec e_z$ implies a prescribed shear $\gamma_{xz}=\Delta L/L$.
Due to symmetry we only consider one half of the domain and prescribe appropriate symmetry boundary conditions at $y=0$.

Two different finite element discretizations are compared. A mesh using prismatic elements aligned with the reinforcement structures, using two elements in each of the stiff layers, and eight respectively four layers in the bulk layers, is used for the full problem. In total, we count 8232 elements in the setup with surface reinforcements, while 10656 elements are counted for the setup with immersed reinforcements. To suppress shear locking, third order displacement elements are used, which results in 341251 respectively 438844 degrees of freedom. An unstructured tetrahedral mesh is designed for the reduced problem, where the stiff layers are approximated as surfaces. The mesh is refined towards the corner points where singularities are expected. For the first setup with surface reinforcements, the mesh consists of 2666 elements, for the second setup 5272 elements are observed. Second order displacement elements (i.e. $k=2$ in \eqref{eq:fe_u0} -- \eqref{eq:fe_m}) are used, resulting in a total of 25240 respectively 45040 degrees of freedom.

\begin{figure}[ht]
  \centering
  \includegraphics[width=0.49\textwidth]{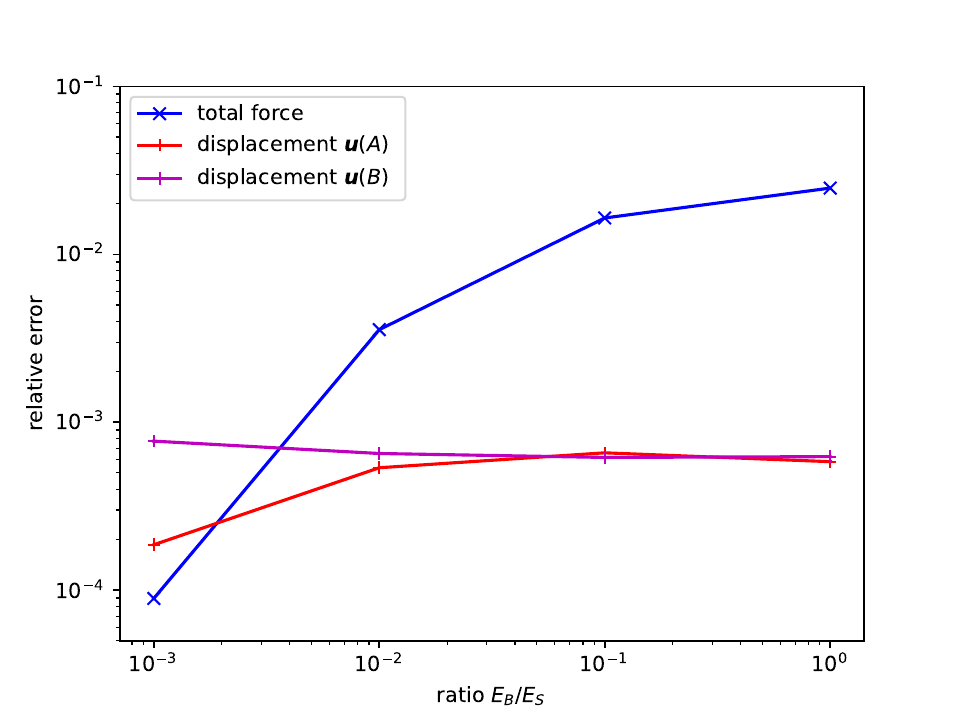}
  \includegraphics[width=0.49\textwidth]{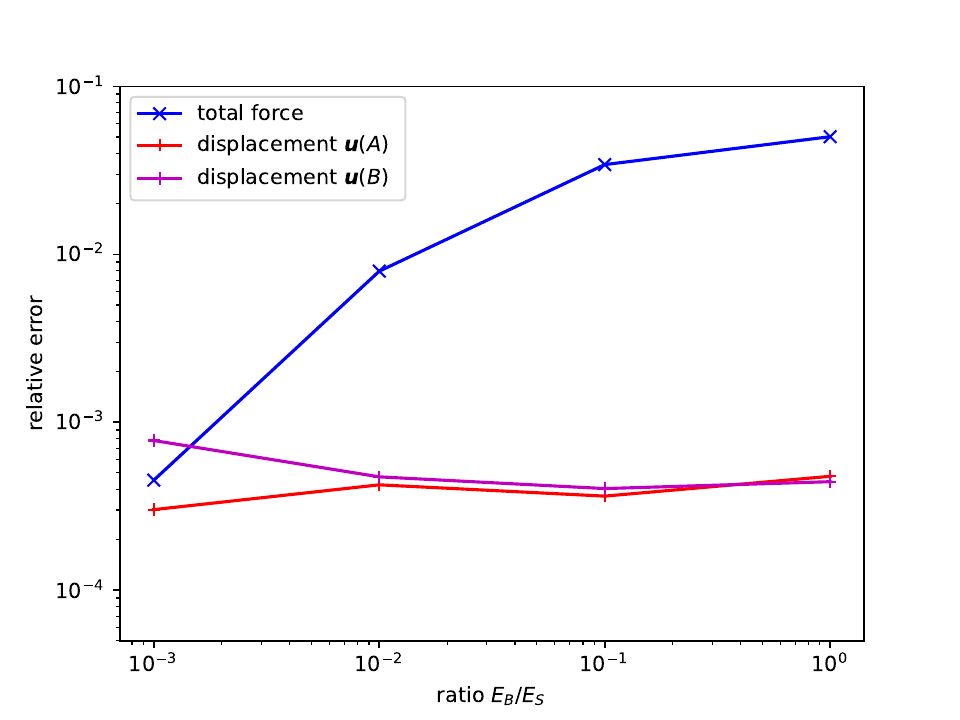}
  \caption{Square laminate plate: load case \emph{extension} for $\Delta L/L = 0.1$, reaction force in longitudinal $x$ direction as well as total displacement $\vec u$ in points $A$ and $B$ are monitored for different stiffness ratios $E_\Bulk/E_\Surf \in [0.001, 1]$. Relative errors of shell discretization with respect to full 3d results are displayed; left: through-the-thickness setup with surface reinforcements; right: through-the-thickness setup with immersed reinforcements.}
  \label{fig:ex1_L1}
\end{figure}

For the first load case (\emph{extension}), a normal strain of $\Delta L/L = 0.1$ is prescribed kinematically. The resulting reaction force is compared for the full three-dimensional and the coupled shell model. Additionally, the displacement of points $A$ and $B$ indicated in Figure~\ref{fig:ex1_geometry} is monitored. Figure~\ref{fig:ex1_L1} presents the relative errors in total force and point displacement for stiffness ratios $E_\Bulk/E_\Surf \in [0.001, 1]$, where in both cases the solution from the fully discretized problem is assumed as reference. We see that for ratios smaller $1:100$, the relative error in reaction force is below $1$\%; for the displacement a relative accuracy of $0.1\%$ is achieved. For higher stiffness ratios, the coupled shell discretization shows larger effective stiffness, and higher computed reaction forces, which is due to the extension of the bulk domain $\Bulk$ up to the shell mid-surfaces. This influence is higher for the submerged reinforcements with a relative error of $5.0\%$ than for the surface reinforcements at $2.5\%$ for $E_\Bulk = E_\Surf$, cf. Figure~\ref{fig:ex1_L1}.

\begin{figure}[ht]
  \centering
  \includegraphics[width=0.49\textwidth]{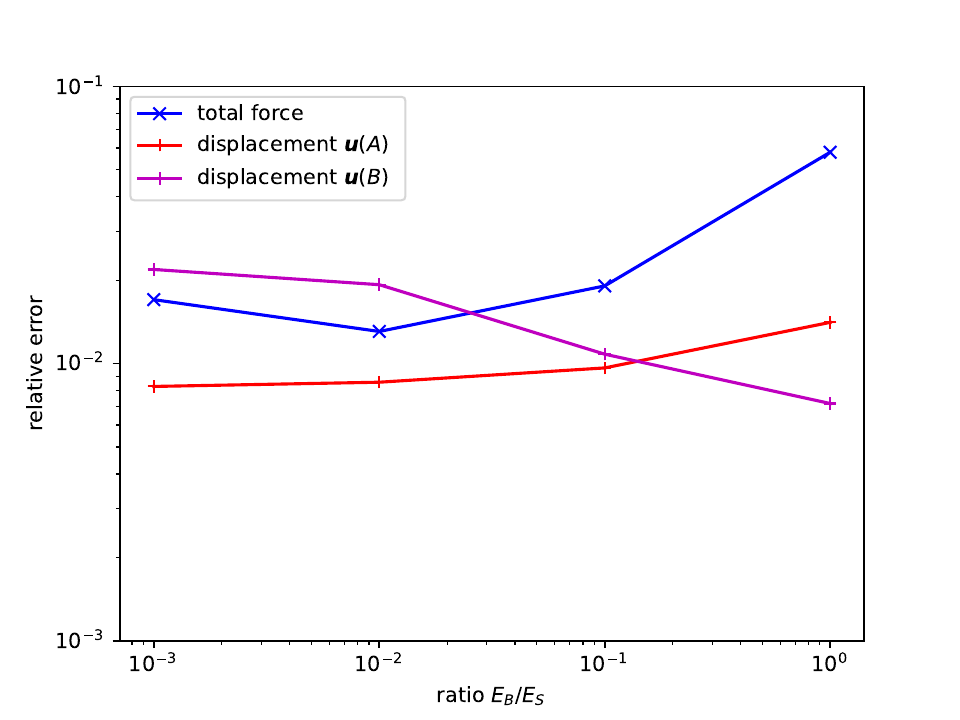}
  \includegraphics[width=0.49\textwidth]{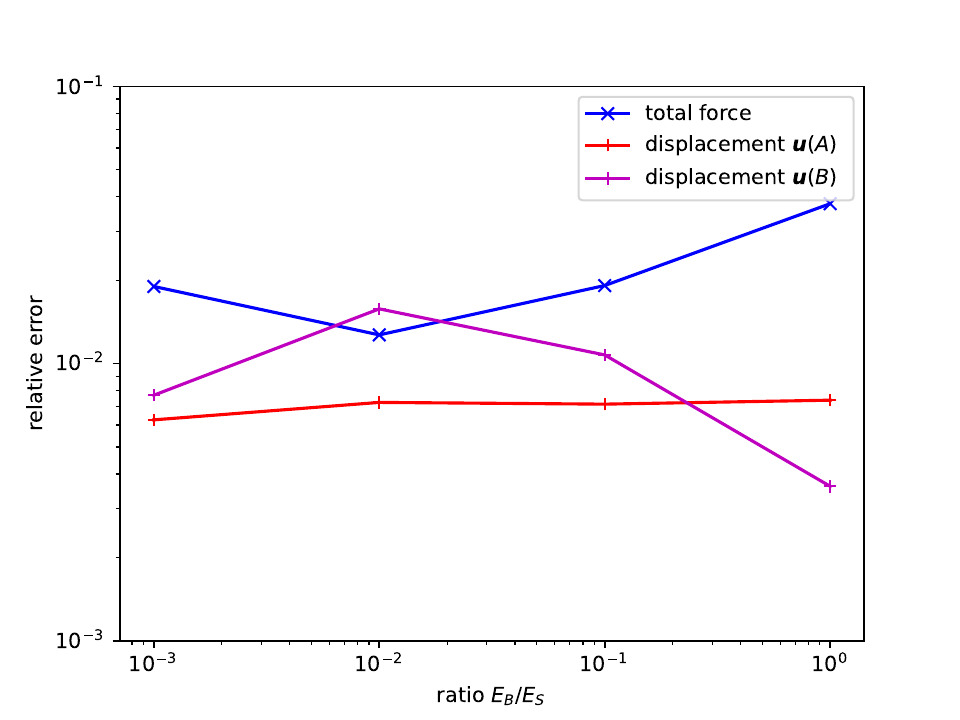}
  \caption{Square laminate plate: load case \emph{simple shear} for $\Delta L/L = 0.1$, reaction force in transverse $z$ direction as well as total displacement $\vec u$ in points $A$ and $B$ are monitored for different stiffness ratios $E_\Bulk/E_\Surf \in [0.001, 1]$. Relative errors of shell discretization with respect to full 3d results are displayed; left: through-the-thickness setup with surface reinforcements; right: through-the-thickness setup with immersed reinforcements.}
  \label{fig:ex1_L2}
\end{figure}

In the second load case (\emph{simple shear}), we choose $\Delta L/L = 0.1$ indicating a shear $\gamma_{xz}$ of $10\%$. The resulting total force in vertical direction is observed for different stiffness ratios $E_\Bulk/E_\Surf \in [0.001, 1]$. Here, the relative error in the reaction force is around $2\%$ for soft bulk material, and increases up to $5.8\%$ for equal stiffness $E_\Bulk = E_\Surf$. The relative error in displacement of points $A$ and $B$ is around $1\%$ for all ratios, see Figure~\ref{fig:ex1_L2}.

\subsection{Half cylinder submerged in soft matrix}
\label{subsec:half_cylinder}

\begin{figure}[!ht]
  \centering
  \begin{tikzpicture}[scale=0.5]
    \def\L{6} 
    \def\H{4} 
    \def\W{20} 
    \def\R{2} 
    \def\opacity{0.5}
    
    \draw[fill=gray!50, opacity=\opacity] (0,0,0) -- ++(\L,0,0) -- ++(0,\H,0) -- ++(-\L,0,0) -- cycle; 
    \draw[fill=gray!50, opacity=\opacity] (0,0,0) -- ++(\L,0,0) -- ++(0,0,\W) -- ++(-\L,0,0) -- cycle; 
    \draw[fill=gray!50, opacity=\opacity] (0,0,0) -- (0,\H,0) -- (0,\H,\W) -- (0,0,\W) -- cycle; 

    \draw[fill=gray!50, opacity=\opacity] (\L,0,0) -- (\L,\H,0) -- (\L,\H,\W) -- (\L,0,\W) -- cycle; 
    \draw[fill=gray!50, opacity=\opacity] (0,\H,0) -- ++(\L,0,0) -- ++(0,0,\W) -- ++(-\L,0,0) -- cycle; 
    \draw[fill=gray!50, opacity=\opacity] (0,0,\W) -- (\L,0,\W) -- (\L,\H,\W) -- (0,\H,\W) -- cycle; 
    \draw[fill=blue!40, opacity=\opacity] (0,0,\W) -- (\L,0,\W) -- (\L,\H,\W) -- (0,\H,\W) -- cycle; 
    \foreach \theta in {180,185,...,355} {
          \draw[fill=blue!40, opacity=\opacity] ({\L/2 + \R*cos(\theta)}, {\H*3/4 + \R*sin(\theta)}, \W) -- ({\L/2 + \R*cos(\theta+5)}, {\H*3/4 + \R*sin(\theta+5)}, \W) -- ({\L/2 + \R*cos(\theta+5)}, {\H*3/4 + \R*sin(\theta+5)}, 0) -- ({\L/2 + \R*cos(\theta)}, {\H*3/4 + \R*sin(\theta)}, 0) -- cycle;
    }    
    \draw[dashed] (0,\H*3/4,0) -- (\L,\H*3/4,0);
    \draw[dashed] (\L/2-\R,\H,0) -- (\L/2-\R,0,0);
    \draw[dashed] (\L/2+\R,\H,0) -- (\L/2+\R,0,0);
    \draw[dashed] (0,\H*3/4,\W) -- (\L,\H*3/4,\W);
    \draw[<->] (\L/2,\H*3/4,0) -- ({\L/2+\R*cos(180+80)},{\H*3/4+\R*sin(180+80)},0) node[midway, right] {$R$};
    \draw[<->] (\L/4,\H*3/4,0) -- ({\L/4},{\H},0) node[midway, right] {$d_1$};
    \draw[<->] (0,\H*9/8,0) -- ({\L/2-\R},{\H*9/8},0) node[midway, above] {$d_2$};

    \draw[<->] (\L,-0.7,\W) -- (0,-0.7,\W) node[midway, below] {$W$};
    \draw[<->] (-0.7,0,\W) -- (-0.7,\H,\W) node[midway, left] {$H$};
    \draw[<->] (\L+0.7,0,\W) -- (\L+0.7,0,0) node[midway, right] {$L$};

    \draw[->] (\L/2,3/4*\H,0) -- (2.2*\L/2,3/4*\H,0) node [anchor=west] {$y$};
    \draw[->] (\L/2,3/4*\H,0) -- (\L/2,1.5*\H,0) node [anchor=west] {$z$};
    \draw[->] (\L/2,3/4*\H,0) -- (\L/2,3/4*\H,1.25*\W) node [anchor=north west]{$x$};
    \end{tikzpicture}
    \caption{Geometry of half cylinder submerged in soft matrix; the tip and the half cylinder constitute a connected shell surface. The shell mid-surface is depicted only.}
    \label{fig:geometry_half_cyl_submerged}
      \end{figure}
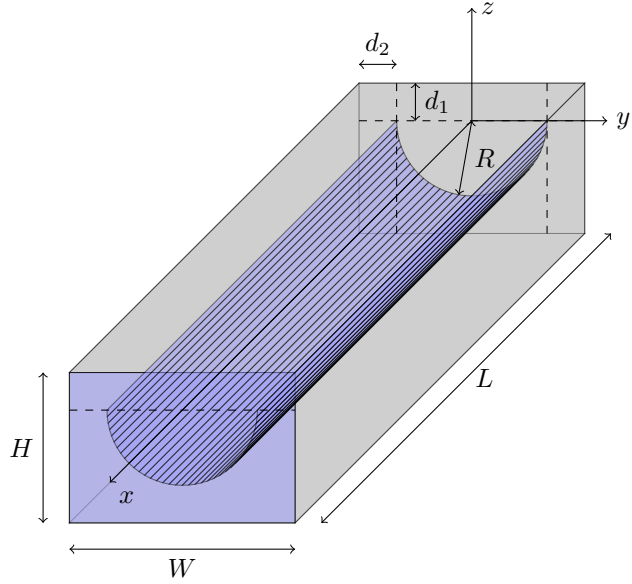
As a second example we consider a half cylinder with radius $R = 2$~mm, thickness $t = 0.1$~mm, and length $L = 20$~mm, which is embedded in a quadrilateral cantilever block of same length $L$, width $W = 6$~mm and height $H = 4$~mm. The half cylinder is located slightly ex-center at distance $d_1 = 1$~mm from the top and $d_2 = 1.1$~mm from the lateral face, cf. Figure~\ref{fig:geometry_half_cyl_submerged}. Additionally, the tip surface is reinforced by a plate of same thickness $t$, which is connected rigidly to the shell. In total, the length of the body in $x$ direction is thus $L + t/2$, similar to the surface reinforcement depicted in Figure~\ref{fig:geometricreduction}. The back face of the cantilever at $x=0$ is fixed rigidly, with the cylindrical shell modeled as clamped. A distributed load is applied in vertical direction to the tip surface of the block, which is denoted as $\Gamma_{\text{tip}}$ in the following. Note that, while in the full 3d simulation this face is located at $x=L+t/2$, in the coupled shell simulation the tip surface coincides with the shell's mid-surface at $x=L$. Depending on the orientation of the load, a buckling point is reached during the application of the load.  To stabilize the numerical procedure in this case, the mean deflection of the tip surface $\bar u_z$ is prescribed via a Lagrangian multiplier $F$, and the resulting force is evaluated, such that
\begin{align}
  \bar u_z &= \frac{1}{|\Gamma_{\text{tip}}|} \int_{\Gamma_{\text{tip}}} \vec u \cdot \vec e_z \, d\Surf, &
  F &= \int_{\Gamma_{\text{tip}}} (\tsr P \cdot \vec N) \cdot \vec e_z\, d\Surf.
\end{align}
The material of the reinforcing half cylinder as well as the tip plate is described by the St.~Venant--Kirchhoff model with $E_\Surf= 70$~GPa and $\nu_\Surf = 0.33$, while the soft matrix material in the cantilever block is given by a compressible Neo-Hooke type with material parameters $E_\Bulk = 10$~MPa and $\nu_\Bulk = 0.2$.

\begin{figure}[ht]
  \centering
  \includegraphics[width=0.48\textwidth]{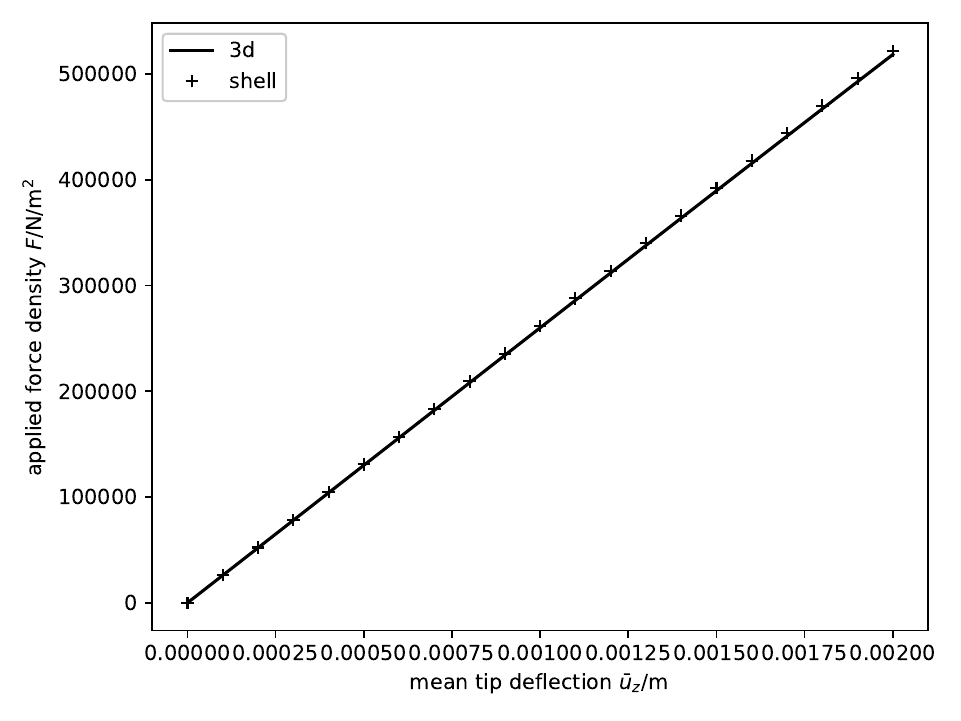}
  \includegraphics[width=0.48\textwidth]{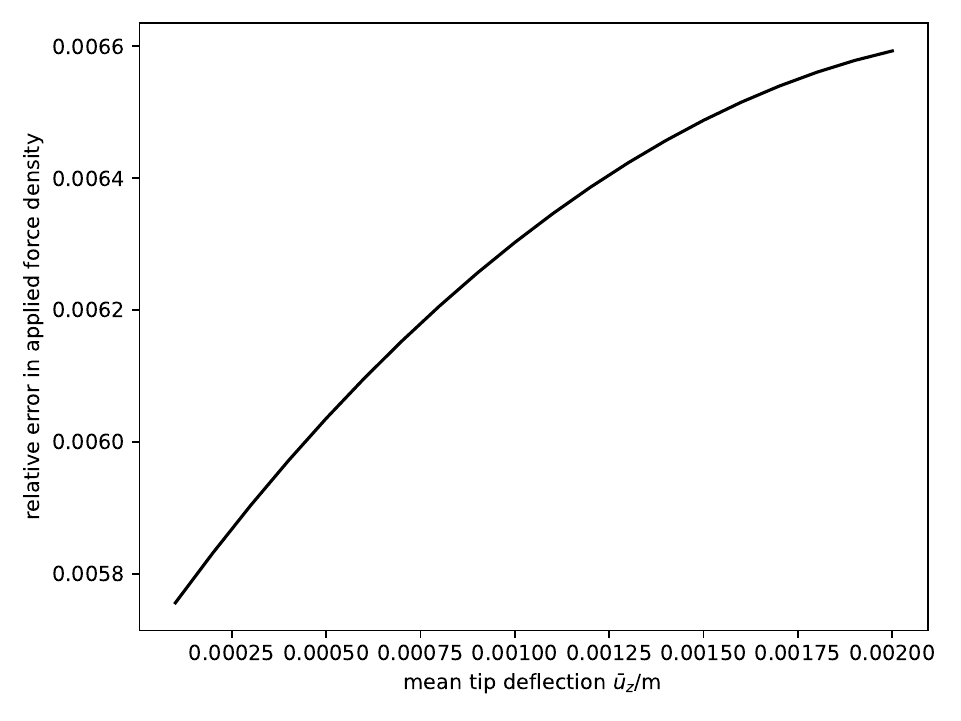}
  \caption{Loadcase without buckling, comparison of tip deflection over load level. Left: mean deflection of tip of bulk material and tip of shell structure. Right: relative error in mean tip deflections.}
  \label{fig:ex2_L1}
\end{figure}

For the first load case, the block face $\Gamma_{\text{tip}}$ at $x=L$ is loaded in negative $z$ direction, the force-displacement relation is nearly linear. The finite element mesh for the full model is an unstructured tetrahedral mesh, that is refined adaptively twice to ensure adequate approximation of the solution. To this end, two steps of iterative refinement according to a stress-based Zienkiewicz-Zhu type error estimate \cite{ZZ1987} is performed, where all elements $T$ with indicator $\eta_T$ above $0.25\eta_{\max}$ are marked for refinement. Using third order basis functions yield $195046$ degrees of freedom for the refined mesh. When using the reduced geometry, only $1742$ elements without adaptive refinement are necessary. Using elements of order two, $19931$ degrees of freedom are used in total, $9203$ of which are inter-element coupling.
We compute the mean displacement of the tip surface, and compare the results for full 3d as well as reduced shell computations at different applied distributed loads. Figure~\ref{fig:ex2_L1} shows an almost linear relation between applied force and mean tip deflection, which is almost perfectly recovered by the reduced shell model. The relative error when taking the full 3d simulation as a reference is below $1\%$.

For the second load case, the tip surface $\Gamma_{\text{tip}}$ is loaded in positive $z$ direction, resulting in buckling at a certain load.
Again, one simulation resolving the original geometry, as well as one with reduced geometry are performed. For the former, third order basis functions are used, while in the latter case, second order elements suffice. In both simulations, adaptive refinement is used to refine the mesh near the points where local buckling occurs. In the coupled shell simulation, independent error estimators for the bulk elements and the shell elements are used. In the two respective cases, all elements $T$ or $T_\Surf$ with indicator $\eta_T \geq 0.5\eta_{\max}$ or $\eta_{T,\Surf} \geq 0.5\eta_{\Surf,\max}$ are marked for refinement. After two steps of adaptive refinement, the mesh for the full model consists of $23388$ tetrahedral elements, and $328330$ degrees of freedom. For the reduced shell model, two steps of iterative refinement lead to only $16660$ elements and $13667$ coupling degrees of freedom are needed.

\begin{figure}[ht]
  \centering
  \includegraphics[width=0.48\textwidth]{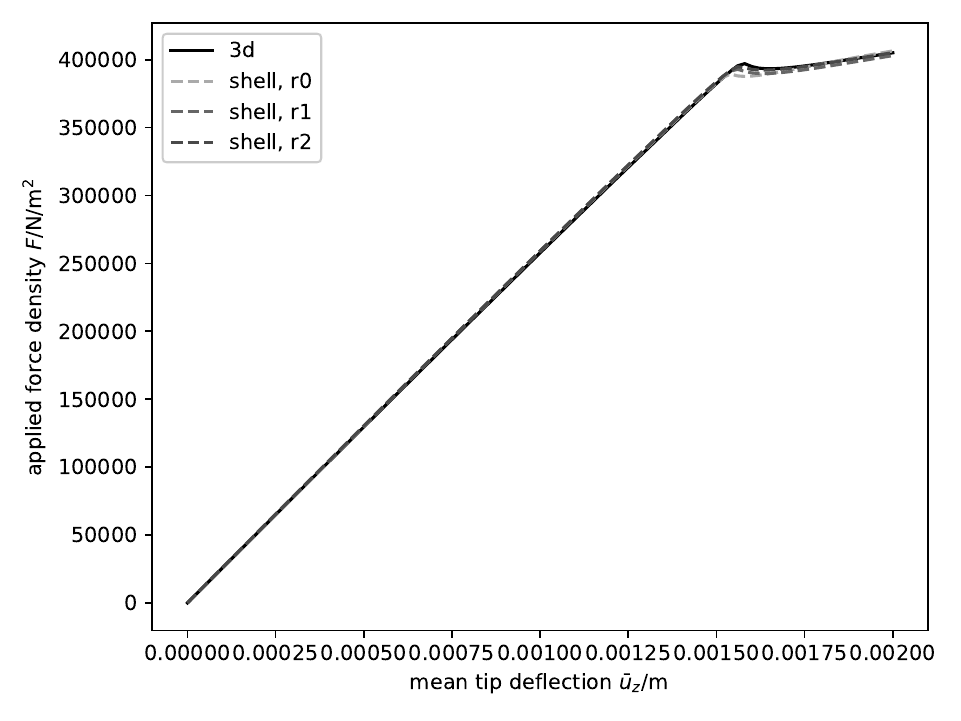}
  \includegraphics[width=0.48\textwidth]{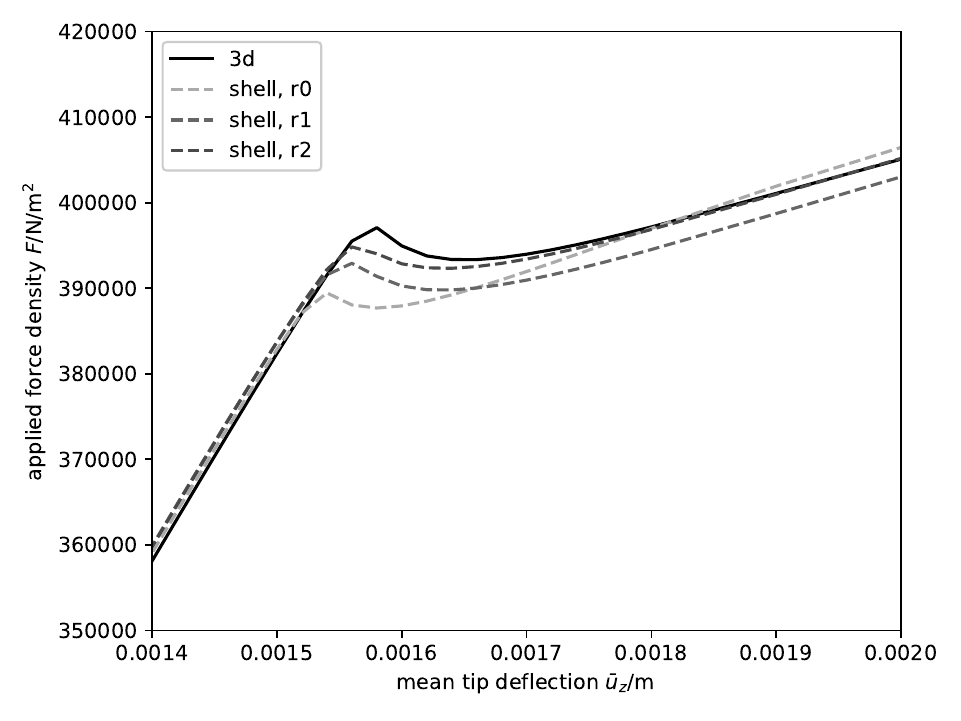}
  \caption{Loadcase with buckling, comparison of tip deflection over load level for different adaptive refinement levels in the shell computation. Right: zoom to the buckling region.}
  \label{fig:ex2_L2_force}
\end{figure}

With sufficient refinement we see almost perfect coincidence in the force level, compare Figure~\ref{fig:ex2_L2_force}. We observe that, once the local refinement is sufficient to resolve the strong deformations in the buckling zone, the full and reduced model yield similar results, the mutual difference way below $1$\%. 
Additionally, we monitor the mean vertical displacement of the block's tip surface in Figure~\ref{fig:ex2_L2_force}. The buckling load, which is computed as $F = 397080$~N at $\bar u_z = 1.58$~mm in the full 3d model, relates well to the load of $F = 394832$~N at $\bar u_z = 1.56$~mm in the reduced shell simulation. In post-buckling regime, the second refinement step in the coupled shell simulation clearly yields more accurate results.

\begin{figure}[ht]
  \centering
  \includegraphics[width=0.4\textwidth]{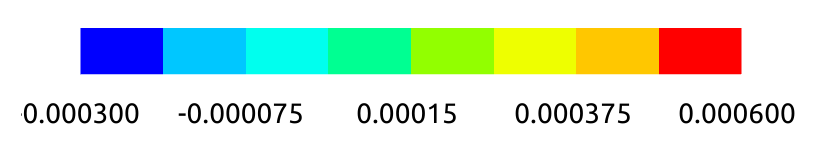}\\
  \includegraphics[width=0.48\textwidth]{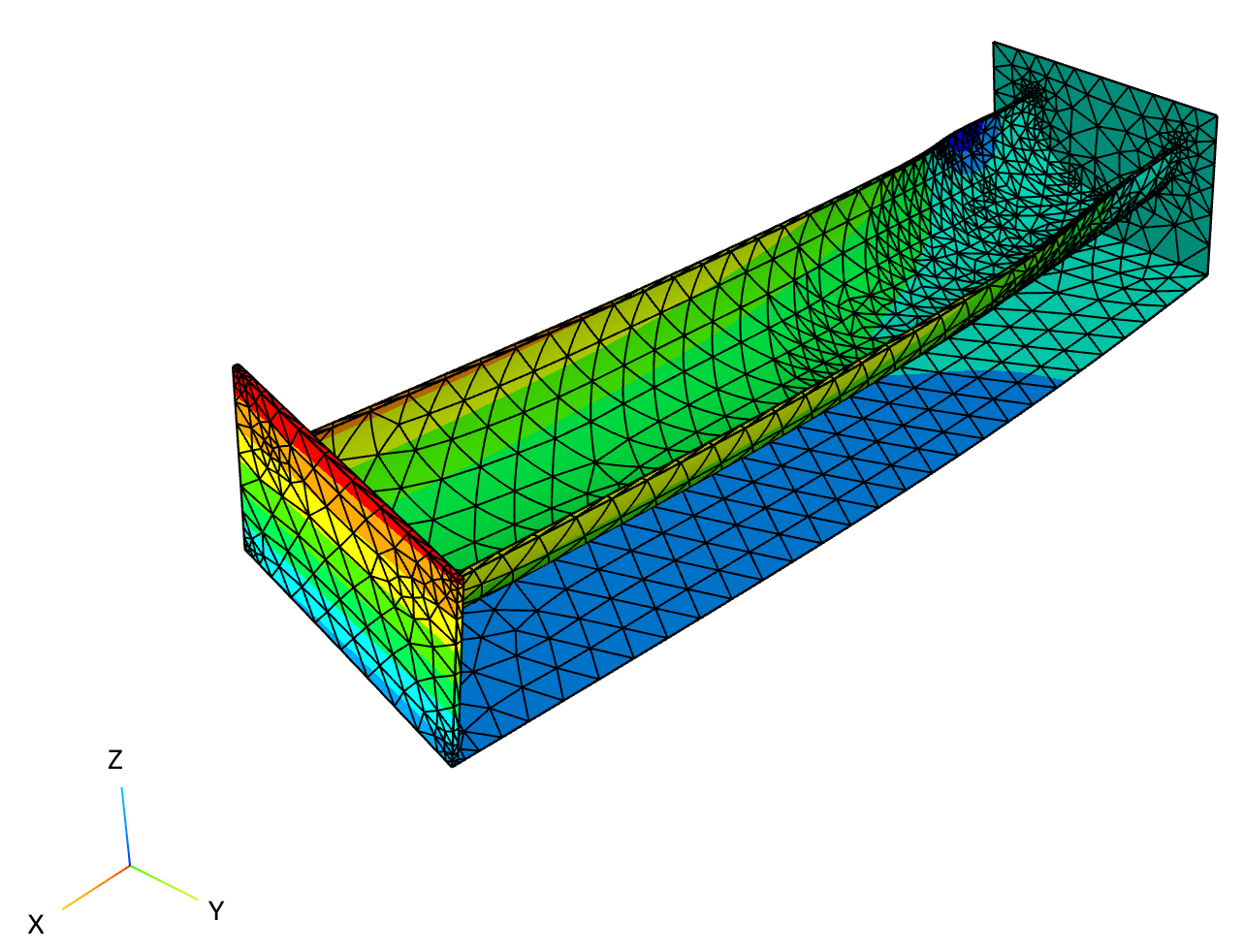}
  \includegraphics[width=0.48\textwidth]{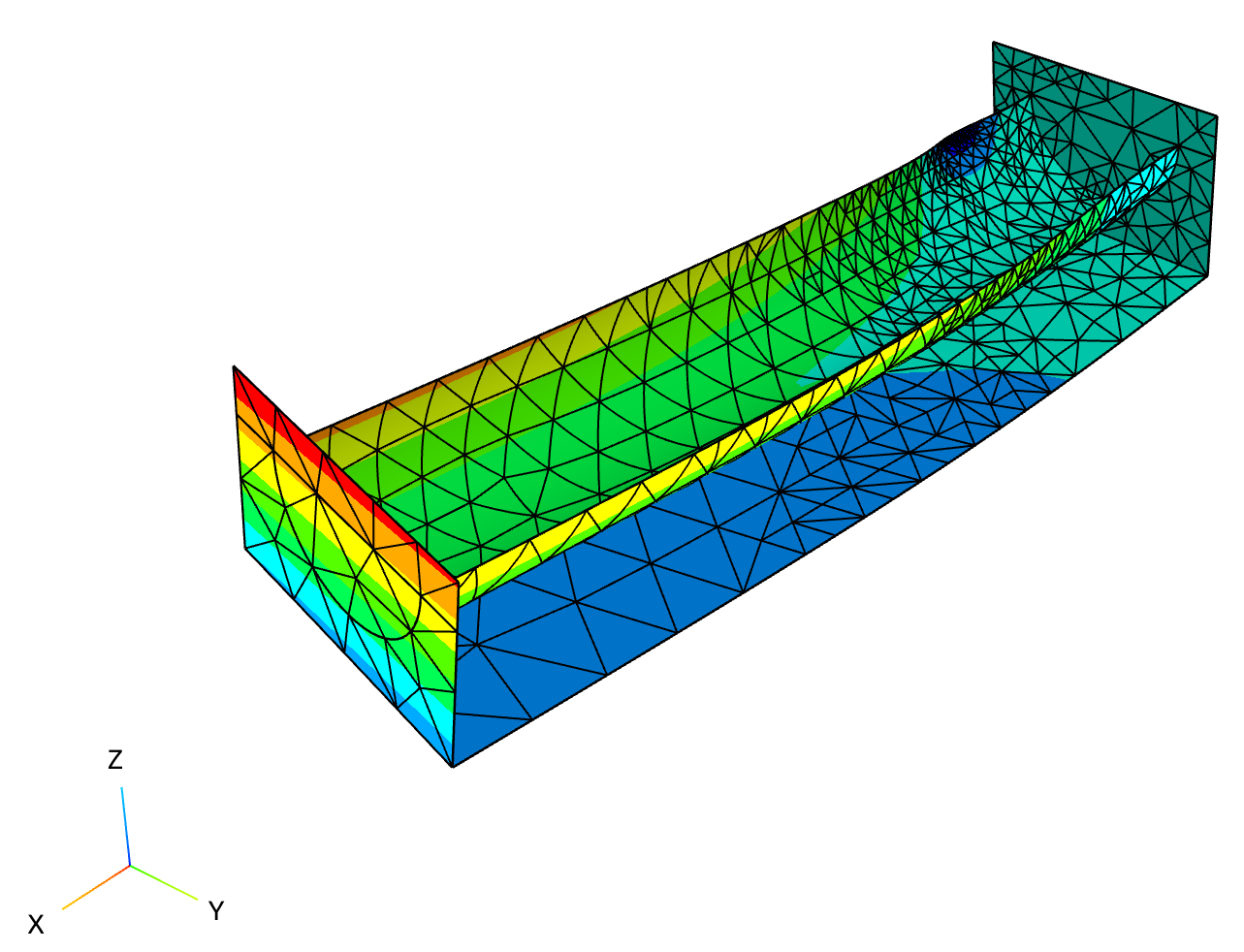}
  \caption{Loadcase with buckling, transverse displacement component $u_y$. Left: full 3d simulation after two adaptive refinement steps; right: coupled shell simulation after two adaptive refinement steps. Only selected surfaces of block and reinforcements are visualized.}
  \label{fig:ex2_L2_uy}
\end{figure}

Figure~\ref{fig:ex2_L2_uy} shows the transverse displacement component $u_y$ evaluated at full deformation of $\bar u_z = 2$~mm. In this plot, only selected surfaces of the bulk and reinforcement materials are visualized. We see that the local buckling behavior is reproduced excellently in the reduced shell model. Adaptive refinement is well suited to detect the buckling region, where the mesh size is reduced accordingly. Further refinement in the coupled shell simulation does not lead to more accurate results, as then the element size locally falls below the shell thickness.

\subsection{Wrinkling problem}
\label{subsec:wrinkling}
In \cite{Knapp:2021}, Knapp and co-workers considered a problem of generating and numerical modeling of microscopic wrinkling surfaces on polydimethylsiloxane (PDMS) elastomers. In their experiment, a block of PDMS is pre-stretched and subsequently part of its surface is plasma-treated. Consequently, a stiffer surface layer evolves, where stiffness and thickness depend on the length of plasma treatment. On release, the bulk material contracts to the original stress-free length, which leads to the formation of wrinkles on the treated surface. Knapp et al.~\cite{Knapp:2021} provide not only measurements for the different wrinkle wavelengths observed, but also data for the thickness and mechanic properties of the surface and bulk material, as well as simulation results. In the current contribution, we set up a similar experiment based on their material data, and test the proposed shell formulation in this demanding setup. Not only do the observed wavelengths make a fine discretization of the treated area necessary; additionally, singularities at the boundaries of the treated surface patches are expected, and the near incompressibility of the bulk material PDMS makes a displacement-pressure formulation in the bulk material necessary.

In our computations, we consider a block of width $W = \SI{20}{\micro\meter}$, height $H = \SI{20}{\micro\meter}$ and length $L = \SI{60}{\micro\meter}$. Centered on the top surface, a patch of width $w = \SI{10}{\micro\meter}$ and length $l = \SI{20}{\micro\meter}$ in undeformed state is chosen for plasma treatment. 
In all computations, the block is assumed to be made from an almost incompressible Neo-Hooke material with $\mu = 0.54$~MPa and $\nu = 0.4995$, which results in a bulk modulus of $\kappa = 560$~MPa. These parameters match the Mooney--Rivlin constants provided in the original reference. In a displacement-pressure formulation, where $\hat {\tsr C} = J^{-2/3} \tsr C$ is the isochoric part of deformation, we consider a mixed potential of the form
\begin{align}
  \psi_\Bulk = \frac{\mu}{2} (\hat{\tsr C} : \tsr I - 3) + p (J-1) - \frac{1}{2\kappa} p^2. 
\end{align}
Three variants of surface treatment are considered, where the surface is exposed to plasma for $1$~min, $2$~min or $4$~min. The respective shell parameters are reported in Table~\ref{tab:wrinkle}. We use unstructured tetrahedral meshes for all computations. In each case, the exposed section of the surface is refined such that the mesh size is approximately one fourth of the expected wave length. This results in finite element meshes with $54621$, $35824$, and $15261$ elements, respectively. The degree of freedom count including the shell discretization is $350572$, $227790$, and $96524$, respectively.

In the simulation, in a first step the homogeneous block is pre-stretched in longitudinal direction; to this end a constant surface load is applied to the two opposite surfaces of the otherwise unconstrained block to achieve a total stretch of $30\%$. In the experiment, the block is plasma-treated in this configuration. In the simulations, the designated part of the top surface is equipped with additional shell degrees of freedom to account for the stiffening. The shell is assumed stress-free in current configuration, which is modeled by including the current membrane strain as eigenstrain; due to the perfect uni-axial deformation of the block, no curvature of the surface is present at this stage.

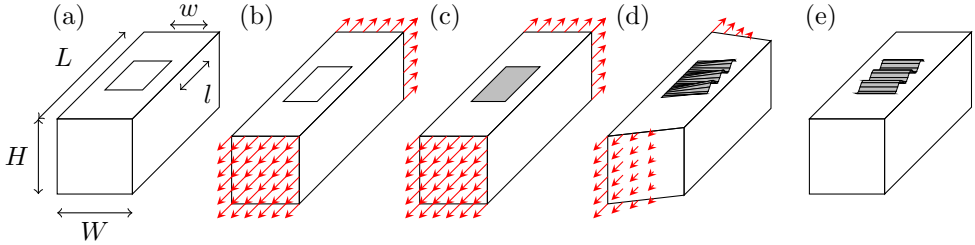
\begin{figure}
  \centering
  \begin{tikzpicture}[scale=0.5]
    \def\L{6} 
    \def\H{2} 
    \def\W{2} 
    \def\l{2} 
    \def\w{1} 
    \def\F{1} 
    \draw (-\W/2,-\H,\L) -- ++(\W,0,0) -- ++(0,\H,0) -- ++(-\W,0,0) -- cycle;
    \draw (\W/2,-\H,0) -- ++(0,0,\L) -- ++(0,\H,0) -- ++(0,0,-\L) -- cycle;
    \draw (-\W/2,0,0) -- ++(0,0,\L) -- ++(\W,0,0) -- ++(0,0,-\L) -- cycle;
    \draw (-\w/2,0,{0.5*(\L-\l)}) -- ++(0,0,\l) -- ++(\w,0,0) -- ++(0,0,-\l) -- cycle;
    \draw[<->] (-\W/2-\W/4,0,0) -- ++(0,0,\L) node [midway,above left] {$L$};
    \draw[<->] (-\W/2-\W/4,-\H,\L) -- ++(0,\H,0) node [midway,left] {$H$};
    \draw[<->] (-\W/2,-\H-\H/4,\L) -- ++(\W,0,0)node [midway,below] {$W$};
    \draw[<->] (-\w/2,0,-\W/4) -- ++(\w,0,0) node [midway,above] {$w$};
    \draw[<->] (\W/2+\W/4,0,{0.5*(\L-\l)}) -- ++(0,0,\l) node [midway,below right] {$l$};
    \node at (-1.5*\W,0.2*\H,0) {(a)};

    \begin{scope}[xshift=5cm]
        \draw (-0.9*\W/2,-0.9*\H,1.2*\L) -- ++(0.9*\W,0,0) -- ++(0,0.9*\H,0) -- ++(-0.9*\W,0,0) -- cycle;
        \draw (0.9*\W/2,-0.9*\H,0) -- ++(0,0,1.2*\L) -- ++(0,0.9*\H,0) -- ++(0,0,-1.2*\L) -- cycle;
        \draw (-0.9*\W/2,0,0) -- ++(0,0,1.2*\L) -- ++(0.9*\W,0,0) -- ++(0,0,-1.2*\L) -- cycle;
        \draw (-0.9*\w/2,0,{1.2*0.5*(\L-\l)}) -- ++(0,0,1.2*\l) -- ++(0.9*\w,0,0) -- ++(0,0,-1.2*\l) -- cycle;
        \foreach \x in {0,...,5}
          \draw [-stealth,color=red] ({-0.9*\W/2*(1-2*\x/5)},0,0) -- ({-0.9*\W*(1-2*\x/5)/2},0,-\F);
        \foreach \x in {0,...,5}
          \draw [-stealth,color=red] ({0.9*\W/2},{-0.9*\H*(\x/5)},0) -- ({0.9*\W/2},{-0.9*\H*(\x/5)},-\F);
        \foreach \x in {0,...,5}
        \foreach \y in {0,...,5}
        \draw [-stealth,color=red] ({-0.9*\W/2*(1-2*\x/5)},{-0.9*\H*(\y/5)},1.2*\L) -- ({-0.9*\W/2*(1-2*\x/5)},{-0.9*\H*(\y/5)},1.2*\L+\F);
    \node at (-1.5*\W,0.2*\H,0) {(b)};
    \end{scope}

    \begin{scope}[xshift=10cm]
      \draw (-0.9*\W/2,-0.9*\H,1.2*\L) -- ++(0.9*\W,0,0) -- ++(0,0.9*\H,0) -- ++(-0.9*\W,0,0) -- cycle;
      \draw (0.9*\W/2,-0.9*\H,0) -- ++(0,0,1.2*\L) -- ++(0,0.9*\H,0) -- ++(0,0,-1.2*\L) -- cycle;
      \draw (-0.9*\W/2,0,0) -- ++(0,0,1.2*\L) -- ++(0.9*\W,0,0) -- ++(0,0,-1.2*\L) -- cycle;
      \draw [fill=gray!50](-0.9*\w/2,0,{1.2*0.5*(\L-\l)}) -- ++(0,0,1.2*\l) -- ++(0.9*\w,0,0) -- ++(0,0,-1.2*\l) -- cycle;
      \foreach \x in {0,...,5}
      \draw [-stealth,color=red] ({-0.9*\W/2*(1-2*\x/5)},0,0) -- ({-0.9*\W*(1-2*\x/5)/2},0,-\F);
    \foreach \x in {0,...,5}
      \draw [-stealth,color=red] ({0.9*\W/2},{-0.9*\H*(\x/5)},0) -- ({0.9*\W/2},{-0.9*\H*(\x/5)},-\F);
    \foreach \x in {0,...,5}
    \foreach \y in {0,...,5}
    \draw [-stealth,color=red] ({-0.9*\W/2*(1-2*\x/5)},{-0.9*\H*(\y/5)},1.2*\L) -- ({-0.9*\W/2*(1-2*\x/5)},{-0.9*\H*(\y/5)},1.2*\L+\F);
    \node at (-1.5*\W,0.2*\H,0) {(c)};

    \end{scope}

    \begin{scope}[xshift=15cm]
      \draw (-0.9*\W/2,-0.9*\H,1.2*\L) -- ++(0.9*\W,0,-0.1*\L) -- ++(0,0.9*\H,0) -- ++(-0.9*\W,0,+0.1*\L) -- cycle;
      \draw (0.9*\W/2,-0.9*\H,0.1*\L) -- ++(0,0,1.*\L) -- ++(0,0.9*\H,0) -- ++(0,0,-1.*\L) -- cycle;
      \draw (-0.9*\W/2,0,0) -- ++(0,0,1.2*\L) -- ++(0.9*\W,0,-0.1*\L) -- ++(0,0,-1.*\L) -- cycle;
      \foreach \x in {0,...,20}
      \draw [fill=gray!50,line width=0](\w/2,{sin(\x/20*1000)*\H/20},{0.5*(\L-\l)+\x/20*\l}) -- (\w/2,{sin((\x+1)/20*1000)*\H/20},{0.5*(\L-\l)+(\x+1)/20*\l}) --
      (-\w/2,0,{0.5*(\L-\l)+(\x+1)/20*\l*1.2}) -- (-\w/2,0,{0.5*(\L-\l)+\x/20*\l*1.2}) -- cycle;
      \foreach \x in {0,...,3}
      \draw [-stealth,color=red] ({-0.9*\W/2*(1-2*\x/5)},0,{0.1*\L*\x/5}) -- ++(0,0,{-\F*(1-\x/5)});
    \foreach \x in {0,...,3}
    \foreach \y in {0,...,5}
    \draw [-stealth,color=red] ({-0.9*\W/2*(1-2*\x/5)},{-0.9*\H*(\y/5)},{1.2*\L-0.1*\L*\x/5}) -- ++(0,0,{\F*(1-\x/5)});
    \node at (-1.5*\W,0.2*\H,0) {(d)};

    \end{scope}

    \begin{scope}[xshift=20cm]
    \draw (-\W/2,-\H,\L) -- ++(\W,0,0) -- ++(0,\H,0) -- ++(-\W,0,0) -- cycle;
    \draw (\W/2,-\H,0) -- ++(0,0,\L) -- ++(0,\H,0) -- ++(0,0,-\L) -- cycle;
    \draw (-\W/2,0,0) -- ++(0,0,\L) -- ++(\W,0,0) -- ++(0,0,-\L) -- cycle;
    \foreach \x in {0,...,20}
      \draw [fill=gray!50,line width=0](\w/2,{sin(\x/20*1000)*\H/20},{0.5*(\L-\l)+\x/20*\l}) -- (\w/2,{sin((\x+1)/20*1000)*\H/20},{0.5*(\L-\l)+(\x+1)/20*\l}) --
      (-\w/2,{sin((\x+1)/20*1000)*\H/20},{0.5*(\L-\l)+(\x+1)/20*\l}) -- (-\w/2,{sin(\x/20*1000)*\H/20},{0.5*(\L-\l)+\x/20*\l}) -- cycle;
    \node at (-1.5*\W,0.2*\H,0) {(e)};

    \end{scope}

  \end{tikzpicture}
  \caption{Wrinkling experiment: idealized setup and simulation steps. (a) dimensions of the block, (b) extension due to axial force, (c) plasma-treatment of surface patch, (d) unsymmetric release of axial force, onset of wrinkling, (e) analyze wrinkling pattern.}
  \label{fig:wrinklingsetup}
\end{figure}

Finally, the axial force acting at the two opposite surfaces is released, which leads to the formation of wrinkles. The force is released in an unsymmetric manner, to introduce a bias in the otherwise symmetric problem, which should alleviate convergence of the numerical solution. All steps of the simulation setup are collected visually in Figure~\ref{fig:wrinklingsetup}.

\begin{figure}
  \centering
  \includegraphics[width=0.48\textwidth]{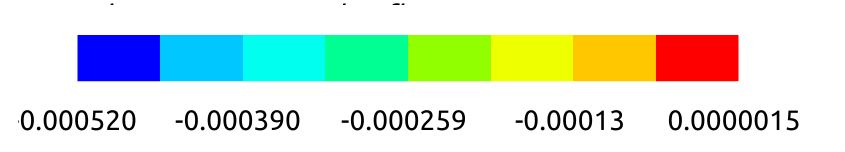}\hfill
  \includegraphics[width=0.48\textwidth]{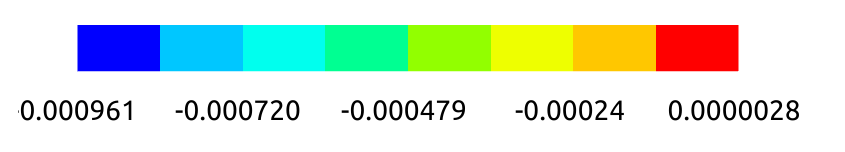}\\
  \includegraphics[width=0.48\textwidth]{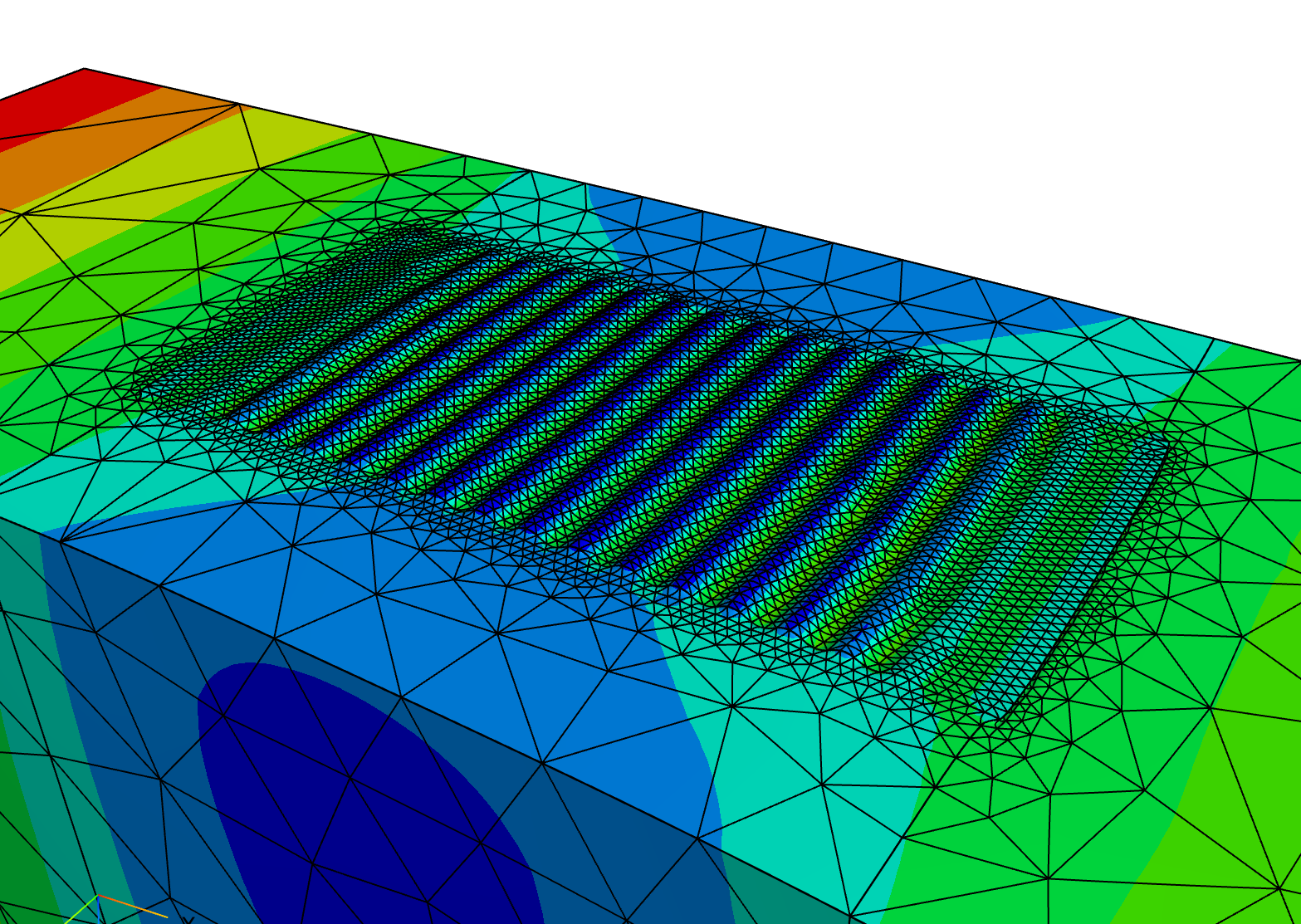}\hfill
  \includegraphics[width=0.48\textwidth]{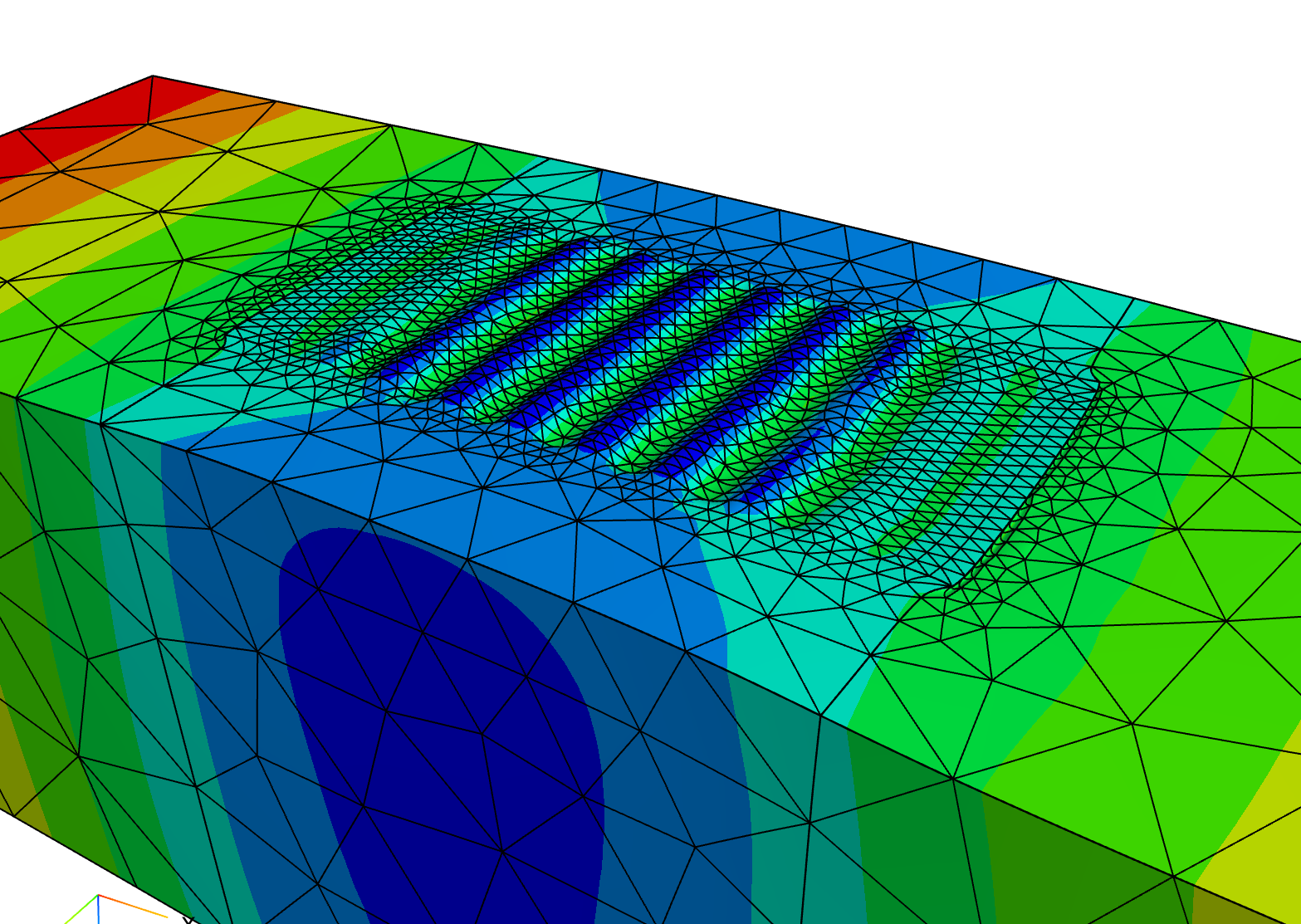}
  \caption{Wrinkling experiment: results of the coupled simulation for exposure times $1$~min (left) and $4$~min (right); visualized is the vertical displacement component $u_z$.}
  \label{fig:wrinklingblock}
\end{figure}

After the axial force is released, the wave pattern of the treated surface is analyzed. Along the center line, wrinkle wave lengths are evaluated. In Table~\ref{tab:wrinkle}, the values obtained in simulations are compared against measurements taken from \cite{Knapp:2021}. Given the uncertainties concerning the mechanical properties of the shell as well the bulk material, as well as the exact geometric and kinematic setup of the silicone block, the wavelengths were reproduced with sufficient accuracy. Figure ~\ref{fig:wrinklingblock} provides a visualization of the computed solution.
\begin{table}
\caption{Wrinkling experiment: shell parameters and wrinkle wavelengths $\lambda_{\exp}$ reported from experiments by \cite{Knapp:2021} as well as wavelengths observed in simulations $\lambda_{\text{sim}}$}
\label{tab:wrinkle}
\begin{center}
  \begin{tabular}{ccccc}
    exposure time & $t_\Surf/\si{\nano\meter}$ & $E_\Surf/\si{MPa}$ & $\lambda_{\exp}/\si{\nano\meter}$ & $\lambda_{\text{sim}}/\si{\nano\meter}$ \\\hline
    $1$~min & $143\pm22$ & $12.1$ & $831 \pm 76$ & 1189\\
    $2$~min & $133\pm5$ & $20.4$ & $1289 \pm 96 $ & 1338\\
    $4$~min & $255\pm33$ & $13.9$ & $1983 \pm 56$ & 2239
    \\\hline
  \end{tabular}
\end{center}
\end{table}

\section*{Summary}
In this paper we presented a novel approach for the coupling of continuum and shell elements in large deformation problems. The proposed method allows for the efficient simulation of structures where thin-walled reinforcements are embedded within or attached to a bulk material. By using a mixed shell formulation, we alleviate the need for higher smoothness in the displacement field, enabling the use of standard nodal elements. Local equilibrium conditions were derived for the coupled problem in the linear case for plates, showing consistency. Several numerical examples were provided to demonstrate the performance and efficiency of the proposed method. Overall, the proposed method provides a powerful tool for the simulation of complex structures involving thin-walled reinforcements, offering significant computational savings while maintaining high accuracy.

\end{document}